\documentclass{article}

\usepackage{amssymb, latexsym}

\usepackage{graphicx}
  
\usepackage{amsmath}
\usepackage{color}

\newtheorem{theorem}{Theorem}[section]

\newtheorem{corollary}[theorem]{Corollary}

\newtheorem{problem}[theorem]{Problem}

\newtheorem{remark}[theorem]{Remark}

\date{ }

\begin{document}
\title{\large\bf Curvature-dimension condition, rigidity theorems and entropy differential inequalities on Riemannian manifolds}
\author{Xiang-Dong Li
\thanks{Research supported by  National Key R$\&$D Program of China (No. 2020YFA0712700), NSFC No. 12571159 
and Key Laboratory RCSDS, CAS, No. 2008DP173182. }}

\maketitle

\begin{minipage}{120mm}

{\bf Abstract} In this paper, we use the information-theoretic approach to study curvature-dimension condition, rigidity 
theorems and entropy differential inequalities on  Riemannian manifolds. We 
prove the equivalence of the ${\rm CD}(K, m)$-condition 
for $K\in \mathbb{R}$ and $m\in [n, \infty]$ and a family of Shannon and R\'enyi entropy differential inequalities along the geodesics on the Wasserstein space over a Riemannian manifold. {The rigidity models of the enhanced entropy differential inequalities are the $K$-Einstein manifolds and the $(K, m)$-Einstein manifolds}. Moreover, we prove the monotonicity and 
rigidity theorem of the $W$-entropy associated with the Shannon entropy and the R\'enyi entropy along the geodesics on the Wasserstein space over Riemannian manifolds with CD$(0, m)$-condition. Comparing with the characterization of the  the CD$(K, m)$ curvature-dimension condition in the framework of the synthetic  geometry  developed by Lott, Sturm and Villani, we provide more simple equivalent characterizations for the CD$(K, m)$-condition, and we provide a characterization of the Einstein and quasi-Einstein manifolds by the enhanced entropy differential equality and the enhanced entropy power differential equality. These are new  in the literature.

\end{minipage}

\medskip

\noindent{\it MSC2010 Classification}: Primary  53C23, 53C21; Secondary 60J60, 60H30.

\medskip

\noindent{\it Keywords}: curvature-dimension condition,  (enhanced) entropy (power) differential inequalities, (quasi-)Einstein metrics, rigidity theorems, $W$-entropy, Lott-Villani-Sturm's synthetic geometry 


%
%

\section{Introduction}
In recent years, there have been extensive works in the study on curvature-dimension condition on Riemannian 
manifolds and metric measure spaces  \cite{AG, BGL, CEMS, EKS, LV, RS, Sturm05, Sturm06a, Sturm06b, V1, V2}. The first link between the $K$-geodesic convexity of the relative entropy functional on the $L^2$-Wasserstein space $(P_2(M), W_2)$ of probability measures on a Riemannian manifold $M$  and the condition on the Ricci curvature bounded below by $K\in \mathbb{R}$ was given by Sturm and von Renesse \cite{RS}. 
In \cite{LV, Sturm06a, Sturm06b}, Lott and Villani, Sturm independently developed a synthetic  comparison geometry on metric measure spaces based on a nice definition of curvature-dimension condition, which extends the notion of the curvature-dimension condition introduced by  Bakry-Emery
 \cite{BE} on Riemannian manifolds with weighted volume measure. By \cite{RS, LV, Sturm06a, V2}, the CD$(K, \infty)$-condition holds on a metric measure space $(X, d, \mu)$ if and only if the Boltzmann (relative) entropy
$${\rm Ent}(\rho)={\rm Ent}(\rho \mu|\mu)=\int_X \rho \log \rho d\mu$$
 is $K$-convex along any geodesic $\{\rho(t)d\mu, t\in [0, 1]\}$ on the Wasserstein space $P_2(M, d)$ over $(X, d, \mu)$
\begin{eqnarray}
{\rm Ent}(\rho(t))\leq (1-t){\rm Ent}(\rho(0))+t{\rm Ent}(\rho(1))-{K\over 2}t(1-t)W_2^2(\rho(0), \rho(1)),
\label{Ent(K)}
\end{eqnarray}
and the CD$(0, N)$-condition holds on a metric measure space $(X, d, \mu)$ if and only if the quantity\footnote{$S_{N'}(\rho)$ is the R\'enyi entropy power $N_{m, p}(\rho)$ for $p=1-{1/N'}$ and $m=N'$. See Section 2. } \begin{eqnarray*}
S_{N'}(\rho)=S_{N'}(\rho \mu|\mu)=-\int_M \rho^{1-1/N'}d\mu
\end{eqnarray*}
 is convex along any geodesic $\{\rho(t)d\mu, t\in [0, 1]\}$ on the Wasserstein space over $(X, d, \mu)$ for all $N'\geq N$,
 \begin{eqnarray}
S_{N'}(\rho(t))\leq (1-t)S_{N'}(\rho(0))+ tS_{N'}(\rho(1)). \label{S(0, N)}
\end{eqnarray}

Recall that the  Wasserstein space $P_2(X, d)$ over $(X, d, \mu)$ is the set of probability measures $\gamma$ on $(X, d)$ with finite second moment $\int_X d^2(o, x)d\gamma(x)<+\infty$ for some (hence for all)  fixed $o\in X$, the Wasserstein distance $W_2(\gamma_0, \gamma_1)$  between two probability measures $\gamma_0, \gamma_1$ on $(X, d)$ is defined by 
$$W_2(\gamma_0, \gamma_1)=\inf\limits_{\pi} \left\{\int_{X\times X} d^2(x, y)d\pi(x, y)\right\}^{1/2}$$
among all the coupling measures $\pi(\cdot, \cdot)$ on $X\times X$ such that $\gamma_0, \gamma_1$ are the marginal measures of $\pi$, i.e., $\gamma_0=\pi(\cdot, X)$ and $\gamma_1=\pi(X, \cdot)$, 
and a geodesic on $P_2(M, d)$ is a curve $\gamma: [0, 1]\rightarrow P_2(M, d)$ such that 
$$W_2(\gamma(s), \gamma(t))=|s-t|W_2(\gamma_0, \gamma_1), \ \ \ \forall s, t\in [0, 1].$$

For $K\neq 0$ and $N< \infty$, the definition formula for the CD$(K, N)$-condition on metric measure spaces is more involved.  In the literature, it was first introduced by Sturm \cite{Sturm06b} using the above mentioned quantity $S_{N'}$ and was further extended by Lott-Villani \cite{LV} using more general functional $U(\rho)=\int_X e(\rho)d\mu$ for continuous convex functions $e$ in the class  $\mathcal{DC}_N$. See Section 4.7 for the precise definition of the class $\mathcal{DC}_N$.

By Definition 1.3 in \cite{Sturm06b}, given two numbers $K, N\in \mathbb{R}$ with $N\in [1, \infty)$, we say that a metric measure
space $(M, d, \mu)$ satisfies the curvature-dimension condition CD$(K, N)$ if and only if for each pair $\nu_0, \nu_1\in P_2(M, d)$ there exist an optimal coupling $q$ of $\nu_0=\rho_0 \mu$ and $\nu_1=\rho_1 \mu$, 
and a geodesic $\Gamma: [0, 1]\rightarrow P_2(M, d)$ connecting $\nu_0$ and $\nu_1$, with
\begin{eqnarray}
S_{N'}(\Gamma(t)|\mu)&\leq &-\int_{M\times M}\left[
\tau^{(1-t)}_{K, N'}(d(x_0, x_1))\rho_0^{-1/N'}(x_0)\right.\nonumber\\
& &\hskip 2cm \left. + \tau^{(t)}
_{K, N'}(d(x_0, x_1))\rho_1^{-1/N'}(x_1)\right]dq(x_0, x_1) 
\label{S(K, N)}
\end{eqnarray}
for all $t\in [0, 1]$ and all $N'\geq N$.  Here, $\tau^{(t)}_{K, N}(\theta)=t^{1/N}\sigma_{K, N-1}^{(t)}(\theta)^{1-1/N}$, with 
\begin{eqnarray*}
\sigma_{K, N}^{(t)}(\theta)={sin( t\theta \sqrt{K/N})\over sin(\theta \sqrt{K/N})}, \ \ t\in [0, 1]
\end{eqnarray*}
when $0<K\theta^2< N\pi^2$, and $\sigma_{0, N}^{(t)}(\theta)=t$ for $K=0$, 
and one replaces the function ${\it sin}$ by the function ${\it sinh}$ in the definition formula of $\sigma_{K, N}^{(t)}(\theta)$ when $K<0$. See also \cite{LV, V2}.

In the CD$(K, N)$-condition on metric measure space, $K$ is regarded as a lower bound of the Ricci curvature and $N\in [1, +\infty)$ is an upper bound of the space dimension. Note that the CD$(K, \infty)$-condition can be viewed as the limiting case of the CD$(K, N)$-condition for  $N\rightarrow \infty$.  For the PDE description of the geodesic $\{\rho(t)d\mu, t\in [0, 1]\}$ on the Wasserstein space over $(X, d, \mu)$, see Section $3$ below.

However,  the above definition formula \eqref{S(K, N)} for the CD$(K, N)$-condition 
looks very complicated. Moreover, it has been an longtime open problem how to use the entropy type formula 
to characterize the Einstein manifolds with constant Ricci curvature and how to extend it to non smooth stetting.

On the other hand, the Shannon and R\'enyi entropies and the corresponding entropy powers have been extensively studied in the information theory. In his 1948 seminal paper \cite{Sh48},  Shannon introduced the differential entropy for an $n$-dimensional random vector $X$ with probability distribution $f(x)dx$ on $\mathbb{R}^n$ by
\begin{eqnarray*}
H(X)=H(f)=-\int_{\mathbb{R}^n} f(x)\log f(x)dx, \label{HX}
\end{eqnarray*}
and the related entropy power by
\begin{eqnarray*}
N(X)=N(f)=e^{{2\over n}H(X)}.\label{NX}
\end{eqnarray*}
In \cite{Sh48}, Shannon discovered the Entropy Power Inequality (EPI): Let $X$ and $Y$ be two independent random vectors with probability density functions 
$f$ and $g$  respectively on $\mathbb{R}^n$, then 
\begin{eqnarray}
N(X+Y)\geq N(X)+N(Y). \label{NXY}
\end{eqnarray}
Equivalently 
\begin{eqnarray*}
N(f*g)\geq N(f)+N(g), \label{Nfg}
\end{eqnarray*}
where $f*g$ denote the convolution of $f$ and $g$, i.e.,  the probability density of the law of $X+Y$.  
See Stam \cite{Stam} and Blachman \cite{Bla} for the complete proofs of the Shannon Entropy Power Inequality.  
Moreover, it is known that the equality holds in $(\ref{NXY})$ if and only if $X$ and $Y$ are normally distributed with proportional covariance matrices. See \cite{Sh48, Stam, Bla, CT}.  

%
%

In \cite{Cost},  Costa proved the Entropy Power Concavity Inequality (EPCI)  for the heat equation on $\mathbb{R}^n$. More precisely, 
let $u$ be a positive solution to the heat equation on $\mathbb{R}^n$
\begin{eqnarray}
\partial_t u=\Delta u. \label{heat1}
\end{eqnarray}
Then the Shannon entropy power $N(u(t))$ is concave in $t\in (0, \infty)$, i.e., 
\begin{eqnarray}
{d^2\over dt^2} N(u(t)) \leq 0. \label{N2}
\end{eqnarray}
Using an argument based on the Blachman-Stam inequality \cite{Bla}, the original proof of Costa's inequalitty $(\ref{N2})$ has been simplified by Dembo et 
al. \cite{CT, Dem1, Dem2}. 
In \cite{V0}, Villani gave a short proof of Costa's  inequality  $(\ref{N2})$ and pointed out that one can extend it to Riemannian manifolds with non-negative Ricci curvature by the Bakry-Emery $\Gamma_2$-calculus.  

The concavity property of the Shannon entropy power for the heat equation $(\ref{heat1})$  has been extended to the R\'enyi entropy power by 
Savar\'e and Toscani \cite{ST}. 
More precisely,  let $u$ be a positive solution to the nonlinear diffusion equation on $\mathbb{R}^n$
\begin{eqnarray}
\partial_t u=\Delta u^p, \label{pme}
\end{eqnarray}
where $p\geq 1-{2\over n}$. Let 
\begin{eqnarray*}
H_p(u)={1\over 1-p}\log \int_M u^pdx
\end{eqnarray*}
be the $p$-th R\'enyi entropy associated with $(\ref{pme})$, and define the R\'enyi entropy power by
\begin{eqnarray}
N_p(u)=\exp\left(\sigma H_p(u)\right),\label{Npu}
\end{eqnarray}
where $\sigma=p-1+{2\over n}$. Then, Savar\'e and Toscani \cite{ST} proved that   the R\'enyi entropy power $N_p(u(t))$ is concave in $t\in (0, \infty)$, i.e., 
\begin{eqnarray}
{d^2\over dt^2} N_p(u(t)) \leq 0. \label{Np}
\end{eqnarray}
Note that when $p\rightarrow 1$, we have $\sigma\rightarrow {2\over n}$,  $H_p(u)\rightarrow H(u)$ and $N_p(u)\rightarrow N(u)$. In view of this, one can recover Costa's concavity result $(\ref{N2})$ from $(\ref{Np})$.

%
%
%
%

In  \cite{LL20a},  Songzi Li and the author of this paper extended the concavity of Shannon entropy power  to complete Riemannian manifolds. 
More precisely, the Shannon entropy power is concave along 
the heat equation associated with the Witten Laplacian on complete Riemannian 
manifolds with the Bakry-Emery curvature-dimension CD$(0, m)$-condition. In particular, the Shannon entropy power  is concave along the heat equation 
$\partial_t u=\Delta u$  on complete Riemannian 
manifolds with non-negative Ricci curvature. Under the curvature-dimension CD$(K, m)$-condition, it is proved in \cite{LL20a} that the rigidity models of the Shannon entropy power are Einstein or quasi Einstein manifolds
with Hessian solitons. 
Moreover, it is also proved in \cite{LL20a}  that  the Shannon entropy power is convex along the conjugate heat equation introduced by G. Perelman for Ricci flow and the rigidity models of the Shannon entropy power are the shrinking Ricci solitons. This gives a new understanding of Perelman's mysterious $W$-entropy functional for the Ricci flow from the information theoretic point of view.

In \cite{Li26}, the author of this paper further proved that the R\'enyi entropy power is concave along 
the  nonlinear diffusion equation $\partial_t u=L u^p$ associated with the Witten Laplacian on compact Riemannian 
manifolds with the Bakry-Emery curvature-dimension CD$(0, m)$-condition. In particular,  for $p\geq 1-{2\over n}$, the R\'enyi entropy power $N_p(u(t))$  is concave along the  nonlinear diffusion equation 
$\partial_t u=\Delta u^p$  on compact Riemannian 
manifolds with non-negative Ricci curvature.   Under the curvature-dimension CD$(K, m)$-condition,  it is proved that the rigidity models of the R\'enyi entropy power are Einstein or quasi Einstein manifolds
with Hessian solitons. 

In the theory of optimal transportation, the heat equation $(\ref{heat1})$  can be realized as the gradient flow of the Boltzmann  entropy on the $L^2$-Wasserstein space 
$P_2^\infty(\mathbb{R}^n)$ of probability densities on $\mathbb{R}^n$ equipped with Otto's infinite dimensional Riemannian metric, 
and the nonlinear diffusion equation $(\ref{pme})$  is the gradient flow of the R\'enyi entropy on the Wasserstein space  over  
$\mathbb{R}^n$.  On the other hand, by the work of McCann \cite{Mc} and Benamou-Brenier \cite{BB},  the geodesic flow  over the $L^2$-Wasserstein space 
$P_2^\infty(\mathbb{R}^n)$  equipped with Otto's Riemannian metric is given by the continuity equation (called also the transport equation) together with the Hamilton-Jacobi equation on $\mathbb{R}^n$. 

The purpose of this paper is use {\it the information-theoretic approach}  to study the 
curvature-dimension condition, rigidity theorems and entropy differential inequalities on complete Riemannian manifolds. We prove the equivalence of the curvature-dimension condition and a family of the Shannon and R\'enyi entropy differential inequalities  along the geodesics flow on the Wasserstein space over Riemannian 
manifolds. We prove the rigidity models of the enhanced entropy differential inequalities are the Einstein manifolds or 
the $(K, m)$-Einstein manifolds.  Moreover, we introduce the $W$-entropy associated with the Shannon entropy and the R\'enyi entropy along 
 the geodesic on the Wasserstein space over Riemannian manifolds and prove its monotonicity and 
rigidity theorem on complete Riemannian manifolds with CD$(0, m)$-condition. Comparing with the characterization of the  the CD$(K, m)$ curvature-dimension condition in the framework of the synthetic  geometry  developed by Lott-Villani \cite{LV}, Sturm \cite{Sturm06a, Sturm06b} and Villani \cite{V2}, we provide more simple equivalent characterizations for the CD$(K, m)$-condition, and we provide a characterization of the Einstein and quasi-Einstein manifolds by the enhanced entropy differential equality and enhanced entropy power differential equality. Our results are new  in the literature.

 \section{Notation and main results}

Let $(M, g)$
be a complete Riemannian manifold, $V\in C^2(M)$ and $d\mu=e^{-V}dv$, where $dv$ is the Riemannian
volume measure on $(M, g)$.  Let $\nabla_\mu^*$ be the $L^2(\mu)$-adjoint of the gradient operator $\nabla$, i.e.,  for any smooth vector field on $M$,
\begin{eqnarray*}
\nabla_\mu^* X=-e^{V}\nabla\cdot (e^{-V}X)=-\nabla \cdot X+\nabla V\cdot  X.
\end{eqnarray*}
The Witten Laplacian acting on smooth functions is defined by
$$L=-\nabla_\mu^* \nabla= \Delta - \nabla V\cdot\nabla.$$
For any $f, \phi\in C^\infty(M)$, we have 
$$\nabla_\mu^* (f\nabla \phi)=-f L\phi-\langle \nabla f, \nabla \phi\rangle.$$
For any $u, v \in C_0^\infty(M)$, the integration by parts formula holds
\begin{eqnarray*}
\int_M \langle \nabla u, \nabla v\rangle d\mu=-\int_M L u vd\mu= - \int_M u L v d\mu.
\end{eqnarray*}
Thus, $L$ is the infinitesimal generator of the Dirichlet form
\begin{eqnarray*}
\mathcal{E}(u, v)=\int_M \langle\nabla u, \nabla v\rangle d\mu, \ \ \ \ u, v\in C_0^\infty(M).
\end{eqnarray*}
By It\^o's theory, the Stratonovich SDE on $M$
\begin{eqnarray*}
dX_t =\sqrt{2}U_t\circ dW_t-\nabla V(X_t)dt,\ \ \ \ \ \nabla_{\circ dX_t} U_t=0,
\end{eqnarray*}
where $U_t$ is the stochastic parallel transport along the trajectory of $X_t$, with initial data $X_0=x$ and $U_0={\rm Id}_{T_xM}$, defines a diffusion process $X_t$ on $M$ with infinitesimal generator $L$.
Moreover,  the  transition probability density of 
the $L$-diffusion process $X_t$ with respect to $\mu$, i.e., the heat kernel  $p_t(x, y)$   of the Witten Laplacian $L$, is the fundamental solution to the heat equation
\begin{eqnarray}
\partial_t u=Lu.  \label{heq}
\end{eqnarray}

In \cite{BE}, Bakry and Emery proved the generalized Bochner formula 
\begin{eqnarray}
L|\nabla u|^2-2\langle \nabla u, \nabla L u\rangle=2|\nabla^2
u|^2+2Ric(L)(\nabla u, \nabla u), \label{BWF}
\end{eqnarray}
where $u\in C_0^\infty(M)$, $\nabla^2 u$ denotes the Hessian of $u$, $|\nabla^2 u|$ is its Hilbert-Schmidt norm, and 
$$Ric(L)= Ric + \nabla^2 V$$
is the so-called (infinite dimensional) Bakry-Emery Ricci curvature associated with the Witten Laplacian $L$.  For $m\in [n, \infty)$,  the $m$-dimensional Bakry-Emery Ricci curvature associated with the Witten Laplacian $L$  is defined by
$$
Ric_{m, n}(L) = Ric + \nabla^2 V - {\nabla V\otimes\nabla V\over m-n}.
$$
In view of this, we have
\begin{eqnarray*}
L|\nabla u|^2-2\langle \nabla u, \nabla L u\rangle \geq {2|Lu|^2\over m}+2Ric_{m, n}(L)(\nabla u, \nabla u).
\end{eqnarray*}
Here we make a convention that $m=n$ if and only if $V$ is a constant. By definition, we have
$$Ric(L)=Ric_{\infty, n}(L).$$
Following \cite{BE}, we say that $(M, g, \mu)$ satisfies the curvature-dimension 
${\rm CD}(K, m)$-condition for a constant $K\in \mathbb{R}$ and $m\in [n, \infty]$ if and only if
$$Ric_{m, n}(L)\geq Kg.$$
Note that, when $m=n$, $V=0$, we have $L=\Delta$ is the usual Laplacian on $(M, g)$, and the ${\rm CD}(K, n)$-condition holds if and only if the Ricci curvature on $(M, g)$ is bounded from below by $K$, i.e., 
$$Ric\geq Kg.$$

Let $P_2(M, \mu)$ (resp. $P_2^\infty(M, \mu)$) be the Wasserstein 
space (resp. the smooth Wasserstein space) of all probability measures $\rho(x)d\mu(x)$ with density function (resp. with smooth density function) $\rho$ on $M$ such that $\int_M d^2(o, x)\rho(x)d\mu(x)<\infty$,  where $d(o, \cdot)$ denotes the distance function from a fixed point $o\in M$.  Similarly to Otto \cite{Ot}, the tangent space $T_{\rho d\mu}P_2^\infty(M, \mu)$ is identified as follows 
\begin{eqnarray*}
T_{\rho d\mu}P_2^\infty(M, \mu)=\{s=\nabla_\mu^*(\rho \nabla \phi): \phi\in C^\infty(M), \ \ \int_M |\nabla \phi|^2\rho d\mu<\infty\},
\end{eqnarray*}
where
$\nabla_\mu^*$ denotes the $L^2$-adjoint of the Riemannian gradient $\nabla$ with respect to the weighted volume measure $d\mu$  on $(M, g)$. 

For $s_i=\nabla_\mu^*(\rho\nabla \phi_i)\in T_{\rho d\mu} P_2^\infty(M, \mu)$, $i=1, 2$, Otto \cite{Ot} introduced the following infinite dimensional Riemannian metric on $P_2^\infty(M, \mu)$ 
\begin{eqnarray*}
\langle \langle s_1, s_2\rangle \rangle:=\int_M \langle\nabla \phi_1, \nabla \phi_2\rangle \rho d\mu,
\end{eqnarray*}
provided that 
\begin{eqnarray*}
\|s_i\|^2:=\int_M |\nabla \phi_i|^2\rho d\mu<\infty, \ \ \ i=1, 2.
\end{eqnarray*}
Let $T_{\rho d\mu}P_2(M, \mu)$ be the completion of $T_{\rho d\mu}P_2^\infty (M, \mu)$ with Otto's Riemannian metric. Then $P_2(M, \mu)$ is an infinite dimensional Riemannian manifold. 

Similarly to  Benamou and Brenier \cite{BB}, we can prove that,  for any given $\mu_i=\rho_i d\mu\in P_2(M, \mu)$, $i=0, 1$, the $L^2$-Wasserstein distance between $\mu_0$ and $\mu_1$ coincides with the geodesic distance between $\mu_0$ and $\mu_1$ on $P_2(M, \mu)$ equipped with Otto's infinite dimensional Riemannian metric, i.e.,  
\begin{eqnarray*}
W_2^2(\mu_0, \mu_1)=\inf\limits\left\{\int_0^1\int_M |\nabla \phi(x, t)|^2\rho(x, t)d\mu(x)dt: \partial_t \rho=\nabla_\mu^*(\rho \nabla \phi), \ \rho(0)=\rho_0, \ \rho(1)=\rho_1\right\}.
\end{eqnarray*}
Given $\mu_0=\rho(\cdot, 0)\mu$, $\mu_1=\rho(\cdot, 1)\mu\in P_2^\infty(M, \mu)$, 
it is known that there is a unique minimizing Wasserstein geodesic $\{\mu(t), t\in [0, 1]\}$  of the form $\mu(t) =(F_t)_*\mu_0$  joining $\mu_0$ and $\mu_1$ in $P_2(M, \mu)$, where $F_t \in {\rm Diff}(M)$  
is given by $F_t(x) = \exp_x(-t \nabla \phi(\cdot, 0))$  for an appropriate Lipschitz function $\phi(\cdot, t)$ (see \cite{Mc}). 
 If the Wasserstein 
geodesic in $P_2(M, \mu)$  belongs entirely to $P_2^\infty(M, \mu)$,  then the geodesic flow $(\rho, \phi)\in TP_2^\infty(M, \mu)$ satisfies the continuity equation  and the Hamilton-Jacobi equation 
\begin{eqnarray}
{\partial_t} \rho-\nabla_\mu^*(\rho \nabla \phi)&=&0,\label{TA}\\
{\partial_t}\phi+{1\over 2}|\nabla \phi|^2&=&0, \label{HJ}
\end{eqnarray}
with the boundary condition $\rho(0)=\rho_0$ and $\rho(1)=\rho_1$.   When  $\rho_0 \in C^\infty(M, \mathbb{R}^+)$ and $\phi_0\in C^\infty(M)$,  defining $\phi(\cdot, t)\in C^\infty(M)$ by the Hopf-Lax solution
 \begin{eqnarray}
\phi(x, t)=\inf\limits_{y\in M}\left(\phi_0(y)+{d^2(x, y)\over 2t}\right),\label{HLS}
 \end{eqnarray} 
 and solving the continuity equation $(\ref{TA})$ by the characteristic method,  it is known that $(\rho, \phi)$ satisfies $(\ref{TA})$ and $(\ref{HJ})$ with $\rho(0)=\rho_0$ and $\phi(0)=\phi_0$. See  \cite{V1}  Sect. 5.4.7. See also \cite{Lo1, Lo2}.
   In view of this, the continuity equation $(\ref{TA})$ 
 and the Hamilton-Jacobi equation $(\ref{HJ})$ describe the geodesic flow on the tangent bundle $TP_2^\infty(M, \mu)$ over the Wasserstein space $P_2(M, \mu)$. 
 Note that the Hamilton-Jacobi equation $(\ref{HJ})$ is also called the eikonnal equation in geometric optics.  In particular, when $m=n$, $V=0$  and $\mu=v$,  we have $\nabla_{v}^*=-\nabla\cdot$, and 
 the geodesic flow on $TP_2(M, v)$  is given by 
 \begin{eqnarray}
\partial_t \rho+\nabla\cdot (\rho \nabla \phi)&=&0, \label{cont1}\\
\partial_t \phi+{1\over 2}|\nabla \phi|^2&=&0. \label{HJ1}
\end{eqnarray}
When $M=\mathbb{R}^n$, let 
\begin{eqnarray}\rho_n(x, t)={e^{-{|x|^2\over 4t^2}}\over (4\pi t^2)^{n/2}}, \ \  
\phi_n(x, t)={|x|^2\over 2t}.
\label{rigiditygeodesic}
\end{eqnarray}
By \cite{LL24},  $(\rho_n, \phi_n)$ is a special solution to the geodesic equations  $(\ref{cont1})$ 
 and $(\ref{HJ1})$ on $P_2(\mathbb{R}^n, dx)$.

Note that our definitions of $P_2(M, \mu)$, $P_2^\infty(M, \mu)$  and the geodesic equations $(\ref{TA})$ 
 and $(\ref{HJ})$ rely on the reference measure $\mu$ on the Riemanian manifold $(M, g)$, while the classical definition of the Wasserstein space $P_2(X,  d)$ relied only on the distance function $d$ on the metric space $X$.

Let $(\rho, \phi)$ be a smooth geodesic  flow on $P_2(M, \mu)$, i.e., a smooth solution to $(\ref{TA})$ and $(\ref{HJ})$. The Shannon entropy associated with the geodesic flow is defined by
\begin{eqnarray*}
H(\rho)=-\int_M \rho\log\rho\hspace{0.2mm} d\mu,
\end{eqnarray*}
 and for $p>1-{1\over m}$, the $p$-th R\'enyi entropy associated with the geodesic flow is defined by
 \begin{eqnarray*}
H_p(\rho)={1\over 1-p}\log \int_M \rho^p\hspace{0.2mm} d\mu.
\end{eqnarray*}
The Shannon entropy power $N_m$ associated with the geodesic flow is defined by
\begin{eqnarray*}
N_m(\rho)=\exp\left({1\over m}H(\rho)\right),
\end{eqnarray*}and the R\'enyi entropy power $N_{m, p}$ associated with the geodesic flow is defined by 
\begin{eqnarray*}
N_{m, p}(u)=\exp\left({1\over m}H_p(u)\right).
\end{eqnarray*}

Let $\theta=W_2(\rho(0)\mu, \rho(1)\mu)$. Define
\begin{eqnarray*}
\sigma_{K, m, \theta}(t):={{\rm sin}( t\theta\sqrt{K/m})\over {\rm sin}(\theta\sqrt{K/m})},
\end{eqnarray*}
when $0<K\theta^2<m\pi^2$, $\sigma_{K, m, \theta(t)}=+\infty$ when $K\theta^2\geq m\pi^2$,  $\sigma_{K, m, \theta}(t)=t$ when $K=0$, and when $K<0$
\begin{eqnarray*}
\sigma_{K, m, \theta}(t)={{\rm sinh}( t\theta \sqrt{K/m})\over {\rm sinh}(\theta \sqrt{K/m})}.
\end{eqnarray*}
Denote
\begin{eqnarray}
\mathcal{N}_{m, K}(t)=\sigma_{K, m, \theta}(1-t)N_m(\rho(0))+ \sigma_{K, m, \theta}(t)N_m(\rho(1)).
\end{eqnarray}
Then 
\begin{eqnarray}
{d^2\mathcal{N}_{m, K}\over dt^2}= -{K\theta^2\over m}\mathcal{N}_{m, K},
 \label{N2m5Kma}
\end{eqnarray}
and 
$$\mathcal{N}_{m, K}(0)=N_m(\rho(0)), \ \ \ \mathcal{N}_{m, K}(1)=N_{m}(\rho(1)).$$

Now we state the main results of this paper.

\begin{theorem}\label{Thm1} Let $(M, g)$ be an $n$-dimensional complete Riemannian manifold, and $K\in \mathbb{R}$.  Let $(\rho, \phi)$ be a smooth solution to the geodesic flow on $TP_2(M, v)$ equipped with 
Otto's infinite dimensional Riemannian metric, i.e., $(\rho, \phi)$ satisfies the continuity equation $(\ref{TA})$ and the Hamilton-Jacobi equation $(\ref{HJ})$ for $\mu=v$. Let $p\geq 1-{1\over n}$  and $d\gamma={\rho^p dv\over \int_M \rho^p dv}$. Then the following conditions are equivalent:
\medskip

(i)  $Ric\geq Kg$,  

(ii) the Shannon entropy differential inequality holds
\begin{eqnarray*}
-H''\geq {1\over n}H'^2+KW_2^2(\rho(0)v, \rho(1)v),\label{nEDIW}
\end{eqnarray*}

(iii) the Shannon entropy power differential inequality holds
\begin{eqnarray*}
{d^2N_n\over dt^2}\leq -{KN_n\over n}W^2_2(\rho(0)v, \rho(1)v),\label{nEPIW}
\end{eqnarray*}

(iv) the R\'enyi entropy differential inequality holds
\begin{eqnarray*}
H_{p}''+{1\over n} H_{p}'^2
\leq -K \int_M |\nabla \phi|^2d\gamma, \label{nREDI}
\end{eqnarray*}

(v) the R\'enyi entropy power differential inequality holds
\begin{eqnarray*}
{d^2\over dt^2}N_{n, p}(\rho(t))
\leq \left(-{K\over n} \int_M |\nabla \phi|^2d\gamma\right)N_{n, p}(\rho(t)),\label{REPIn}
\end{eqnarray*}

(vi) the Shannon entropy power $N_n(\rho(t))$ satisfies  

\begin{eqnarray*}
N_n(\rho(t))\geq \sigma_{K, n, \theta}(1-t)N_n(\rho(0))+ \sigma_{K, n, \theta}(t)N_n(\rho(1)),\label{NnKk}
\end{eqnarray*}

(vii) For all $N\geq n$, it holds
 
 \begin{eqnarray*}
 {d^2\over dt^2}S_N(\rho)+
{N-n\over n}\left[S_{N}(\rho)\right]^{-1}\left({d\over dt}S_N(\rho)\right)^2
\geq {K\over N}\int_M |\nabla \phi|^2 \rho^{1-1/N}dv,\label{SNn4a}
\end{eqnarray*}
where $S_N(\rho)=-\int_M \rho^{1-1/N}dv$, 
in particular, it holds

 \begin{eqnarray*}
 {d^2\over dt^2}S_N(\rho)\geq {K\over N}\int_M |\nabla \phi|^2 \rho^{1-1/N}dv,\label{SNn4b}
\end{eqnarray*}

(viii) Sturm's definition inequality $(\ref{S(K, N)})$ for CD$(K, N)$ holds for all $N'\geq N=n$.
%
%
\end{theorem}

More generally, we have the following 

\begin{theorem}\label{Thm2} Let $(M, g)$ be  an $n$-dimensional complete Riemannian manifold with a weighted volume measure $d\mu=e^{-V}dv$, where $V\in C^2(M, \mathbb{R})$. Let $m\in [n, \infty)$ and $K\in \mathbb{R}$ be two constants.
Let $(\rho, \phi)$ be a smooth solution to the geodesic flow on $TP_2(M, \mu)$ equipped with 
Otto's infinite dimensional Riemannian metric, i.e., $(\rho, \phi)$ satisfies the continuity equation $(\ref{TA})$ and the Hamilton-Jacobi equation $(\ref{HJ})$. Let  $p\geq 1-{1\over m}$   and $d\gamma={\rho^p d\mu\over \int_M \rho^p d\mu}$. Then the following conditions are equivalent:

\medskip

(i)  $Ric_{m, n}(L)\geq Kg$,  i.e., the CD$(K, m)$-condition holds,

(ii) the Shannon entropy differential inequality holds
\begin{eqnarray}
-H''\geq {1\over m}H'^2+KW_2^2(\rho(0)\mu, \rho(1)\mu),\label{EDIW}
\end{eqnarray}

(iii) the Shannon entropy power differential inequality holds
\begin{eqnarray}
{d^2N_m\over dt^2}&\leq& -{KN_m\over m}W^2_2(\rho(0)\mu, \rho(1)\mu), \label{EPIW}
\end{eqnarray}

(iv) the R\'enyi entropy differential inequality holds
\begin{eqnarray}
H_p''+{1\over m} H_p'^2
\leq -K \int_M |\nabla \phi|^2d\gamma, \label{REDI}
\end{eqnarray}

(v) the R\'enyi entropy power differential inequality holds
\begin{eqnarray}
{d^2\over dt^2}N_{m, p}(\rho(t))
\leq \left(-{K\over m} \int_M |\nabla \phi|^2d\gamma\right)N_{m, p}(\rho(t)),\label{REPI}
\end{eqnarray}

(vi) the Shannon entropy power $N_m(\rho(t)$ satisfies  

\begin{eqnarray}
N_m(\rho(t))\geq \sigma_{K, m, \theta}(1-t)N_m(\rho(0))+ \sigma_{K, m, \theta}(t)N_m(\rho(1)),\label{NmKk}
\end{eqnarray}

(vii) For all $N\geq m$, it holds
 
 \begin{eqnarray}
 {d^2\over dt^2}S_N(\rho)+
{N-m\over m}\left[S_{N}(\rho)\right]^{-1}\left({d\over dt}S_N(\rho)\right)^2
\geq {K\over N}\int_M |\nabla \phi|^2 \rho^{1-1/N}d\mu,\label{SNm4a}
\end{eqnarray}
where $S_N(\rho)=-\int_M \rho^{1-1/N}d\mu$, 
in particular, it holds

 \begin{eqnarray}
 {d^2\over dt^2}S_N(\rho)\geq {K\over N}\int_M |\nabla \phi|^2 \rho^{1-1/N}d\mu,\label{SNm4b}
\end{eqnarray}

(viii) Sturm's definition inequality   $(\ref{S(K, N)})$ holds for all  $N'\geq N=m$. 
\end{theorem}

\begin{theorem}\label{rigidity3} Under the same notations as in Theorem \ref{Thm2}, we have
\begin{eqnarray}
H_p''+{1\over m} H_p'^2
&=& \left(  p-1+{1\over m} \right) \left[\left(  \int_M L \phi d\gamma \right)^2-   \int_M  |L \phi|^2 d\gamma \right]-\int_M Ric_{m, n}(L)(\nabla \phi, \nabla \phi) d\gamma\nonumber\\
& &\hskip0.5cm -\int_M \left[ {m-n\over mn}\left(L \phi+{m\over m-n}\nabla V\cdot\nabla \phi\right)^2 + \left\|\nabla^2 \phi-{\Delta \phi\over n} g\right\|^2_{\rm HS} \right]d\gamma,\nonumber\\ 
&& \label{Hpm2}
\end{eqnarray}  
and
\begin{eqnarray}
{d^2\over dt^2}N_{m, p}(\rho(t))
&=&-\left(  p-1+{1\over m} \right) { N_{m, p}(\rho(t))\over m}{\rm Var}_\gamma(L \phi)-{N_{m, p}(\rho(t))
\over m}\int_M Ric_{m, n}(L)(\nabla \phi, \nabla \phi) d\gamma\nonumber \\
& &-{N_{m, p}(\rho(t))\over m} \int_M \left[ {m-n\over mn}\left(L \phi+{m\over m-n}\nabla V\cdot\nabla \phi\right)^2 + \left\|\nabla^2 \phi-{\Delta \phi\over n} g\right\|^2_{\rm HS} \right]d\gamma,\nonumber\\
& & \label{Np2}
\end{eqnarray}
where
\begin{eqnarray*}
{\rm Var}_\gamma(L \phi)=\int_M |L \phi-\gamma(L \phi)|^2d\gamma.
\end{eqnarray*}
{In particular, under the CD$(K, m)$-condition, i.e.,  $Ric_{m, n}(L)\geq Kg$ with $K\in C(M, \mathbb{R})$,  
and $p\geq 1-{1\over m}$, we have the enhanced entropy differential inequality
\begin{eqnarray}
& &H_p''+{1\over m} H_p'^2+ \left(  p-1+{1\over m} \right){\rm Var}_\gamma(L\phi)\nonumber\\ 
& &\hskip0.5 cm+ \int_M \left[{m-n\over mn}\left(L \phi+{m\over m-n}\nabla V\cdot\nabla \phi\right)^2 + 
\left\|\nabla^2 \phi-{\Delta \phi\over n} g\right\|^2_{\rm HS} \right] d\gamma\nonumber\\ 
& & \hskip3cm \leq -\int_M K |\nabla \phi|^2d\gamma, \label{Hpm2K}
\end{eqnarray}  
and the enhanced entropy power differential inequality
\begin{eqnarray}
& &{d^2\over dt^2}N_{m, p}(\rho(t))+\left(  p-1+{1\over m} \right) { N_{m, p}(\rho(t))\over m}{\rm Var}_\gamma(L \phi)\nonumber\\
& &\hskip1cm +{N_{m, p}(\rho(t))\over m} \int_M \left[ {m-n\over mn}\left(L \phi+{m\over m-n}\nabla V\cdot\nabla \phi\right)^2 + \left\|\nabla^2 \phi-{\Delta \phi\over n} g\right\|^2_{\rm HS} \right]d\gamma\nonumber\\
& &\hskip4cm \leq -{N_{m, p}(\rho(t))\over m}\int_M K|\nabla \phi|^2 d\gamma.
 \label{Np2K}
\end{eqnarray}
In particular, we have
\begin{eqnarray}
{d^2\over dt^2} N_{m, p}(\rho(t))\leq -{N_{m, p}(\rho(t)) \over m} \int_M K |\nabla \phi|^2 d\gamma. \label{KNp1}
\end{eqnarray}
Moreover, the equality in the enhanced entropy differential inequality \eqref{Hpm2K} or the enhanced entropy power differential inequality $(\ref{Np2K})$ holds for all geodesics $(\rho, \phi)$ at some $t$ if and only if $(M, g, V)$ is $(K, m)$-Einstein, i.e., 
$$Ric_{m, n}(L)=Kg.$$ 
Furthermore, in the case $Ric_{m, n}(L)=Kg$, \eqref{Hpm2K}  or \eqref{KNp1} becomes an equality if and only if $\phi$ satisfies the Hessian soliton equation
\begin{eqnarray*}
L\phi=I_p=\int_M L\phi d\gamma, \ \ L \phi+{m\over m-n}\nabla V\cdot\nabla \phi=0, \ \ \nabla^2 \phi={I_p\over m}g.
\end{eqnarray*}
}
\end{theorem}

In particular, for $m=n$, $V=0$, $\mu=v$, $L=\Delta$, we have  

\begin{theorem} \label{rigidity5}Let $N_{n, p}$ be the $p$-th entropy power associated with the R\'enyi entropy for the geodesic flow on $TP_2(M, v)$, and $d\gamma={\rho^p dv\over \int_M \rho^p dv}$. Then
\begin{eqnarray*}
{d^2\over dt^2}N_{n, p}(\rho(t))
&=&- \left(  p-1+{1\over n} \right)  {N_{n, p}(\rho(t))\over n}{\rm Var}_\gamma(\Delta \phi)- {N_{n, p}(\rho(t))\over n}\int_M Ric(\nabla \phi, \nabla \phi) d\gamma\\
& &\hskip2cm - {N_{n, p}(\rho(t))\over n} \int_M \left\|\nabla^2 \phi-{\Delta \phi\over n} g\right\|^2_{\rm HS}d\gamma,
\end{eqnarray*}
where
\begin{eqnarray*}
{\rm Var}_\gamma(\Delta \phi)=\int_M \left|\Delta \phi-\int_M \Delta \phi d\gamma\right|^2d\gamma.
\end{eqnarray*}
In particular, when $Ric\geq Kg$ with $K\in C(M, \mathbb{R})$ and $p\geq 1-{1\over n}$, we have
\begin{eqnarray}
{d^2\over dt^2}N_{n, p}(\rho(t))\leq - {N_{n, p}(\rho(t))\over n}\int_M K |\nabla \phi|^2d\gamma, \label{NpK}
\end{eqnarray}
and
\begin{eqnarray}
 {d \over dt}W_{n, p}(\rho, t)\leq -t\int_M K|\nabla\phi|^2d\gamma-{t\over n} \left|I_p-{n\over t}\right|^2,  \label{NIW-2n}
\end{eqnarray}
where $$I_p=\int_M \Delta \phi d\gamma.$$

Moreover,  under the assumption $Ric\geq Kg$ for $K\in C(M, \mathbb{R})$ and $p\geq 1-{1\over n}$, we have the enhanced entropy power differential inequality
\begin{eqnarray}
& &{d^2\over dt^2}N_{n, p}(\rho(t))+ \left(  p-1+{1\over n} \right)  {N_{n, p}(\rho(t))\over n}{\rm Var}_\gamma(\Delta \phi)+ {N_{n, p}(\rho(t))\over n} \int_M \left\|\nabla^2 \phi-{\Delta \phi\over n} g\right\|^2_{\rm HS}d\gamma\nonumber\\
& &\hskip4cm \leq - {KN_{n, p}(\rho(t))\over n}\int_M |\nabla \phi|^2d\gamma. \label{NIWK}
\end{eqnarray}
The equality in the enhanced entropy power differential inequality $(\ref{NIWK})$ holds for all geodesics $(\rho, \phi)$ at some $t$ if and only if $(M, g)$ is $K$-Einstein in the sense that
$$Ric=Kg.$$  
Furthermore, under the condition $Ric=Kg$, \eqref{NpK} or \eqref{NIW-2n} become an equality if and only if $\phi$ satisfies the Hessian soliton equation
\begin{eqnarray*}
\Delta \phi=I_p:=\int_M \Delta \phi d\gamma, \ \ \nabla^2 \phi={I_p\over n}g.
\end{eqnarray*}

\end{theorem}

Finally, we have the following 

\begin{theorem}\label{W-entropy}  
Let $(M, g)$ be a complete Riemannian manifold and $V\in C^2(M)$ such that $Ric(L)=Ric+\nabla^2V$ is uniformly bounded on $M$. Let $(\rho(t), \phi(t))$ be a smooth Benamou-Brenier geodesic on $P_2(M, \mu)$ satisfying the  growth condition as required in Theorem \ref{AAAA} below. Let $p\geq 1-{1\over m}$   and $d\gamma={\rho^p d\mu\over \int_M \rho^p d\mu}$. Then
\begin{eqnarray*}
\ \ \ \ & & {d^2\over dt^2}H_p(\rho(t))+{2\over t} {d\over dt}H_p(\rho(t))-{m\over t^2} \\
&=&-t(p-1){\rm Var}_\gamma(L \phi)-t\int_M Ric_{m, n}(L)(\nabla \phi, \nabla \phi)d\gamma\\
& & -t \int_M \left(\left|\nabla^2 \phi-{g\over t}\right|^2+{1\over m-n}\int_M \left|\nabla V\cdot \nabla \phi+{m-n\over t}\right|^2\right) d\gamma,
\end{eqnarray*}
where 
\begin{eqnarray*}
{\rm Var}_\gamma(L \phi)=\int_M |L \phi-\gamma(L \phi)|^2d\gamma.
\end{eqnarray*}
Let
\begin{eqnarray*}
H_{m, p}(\rho(t)|\rho_m(t))=H_p(\rho(t))-H_p(\rho_m(t)),
\end{eqnarray*}
and define the $W$-entropy by the Boltzmann formula
\begin{eqnarray*}
W_{m, p}(\rho, t)={d\over dt}(tH_{m, p}(\rho(t)|\rho_m(t)).
\end{eqnarray*}
Then
 \begin{eqnarray*}
{1\over t}{d\over dt}W_{m, p}(\rho, t)
&=&-(p-1){\rm Var}_\gamma(L \phi)- \int_M  Ric_{m, n}(L)(\nabla \phi, \nabla \phi) d\gamma\nonumber\\
& & -\int_M \left[\left\|\nabla^2 \phi-{g\over t}\right\|^2 +{1\over m-n}\int_M \left(\nabla V\cdot\nabla \phi+{m-n\over t}\right)^2\right]d\gamma.\nonumber\\
& &\label{derivativeofW}
\end{eqnarray*}
In particular, if $Ric_{m, n}(L)\geq 0$, then $tH_{m, p}(\rho(t))$ is concave in $t$ along the geodesic in $P^\infty(M, \mu)$, and $W_{m, p}(\rho(t))$ is non-increasing in $t$ 
along the geodesic in $P^\infty(M, \mu)$. 

\medskip
Moreover, assume that  
$Ric_{m, n}(L)\geq 0$, $p\geq 1-{1\over m}$, and ${d\over dt}W_{m, p}(\rho, t)=0$ at some $t=\tau>0$. Then $M$ is isometric to $\mathbb{R}^n$, $m=n$, $V$ is a constant and $(\rho, \phi)=(\rho_n, \phi_n)$ as given by \eqref{rigiditygeodesic}.

\end{theorem}

For the significant role of the $W$-entropy in the proof of the non-local collapsing theorem of Hamilton's Ricci flow \cite{H82, H95} 
for  the Poincar\'e conjecture and Thurston's geometrization conjecture, see Perelman \cite{P1}. See also \cite{CZ, KL, MT, CLN, CCLLN}.  See \cite{N1, N2} for the $W$-entropy formula for the heat equation $\partial_tu=\Delta u$ on Riemannian manifolds with non-negative Ricci curvature.  
See \cite{Li12, LL15, LL18a, LL18, LL24} for the study of monotonicity and rigidity theorems of the  $W$-entropy for the heat equation of the Witten Laplacian $\partial_tu=L u$ on Riemannian manifolds with CD$(K, m)$-condition or $(K, m)$-super Ricci flows and 
its extension to the Langevin deformation on the Wasserstein space over Riemannian mnaifolds with CD$(0, m)$-condition.

Note that, in the work of Sturm and von Renesse \cite{RS}, Lott-Villani \cite{LV} and Sturm \cite{Sturm06a, Sturm06b}, the $K$-displacement convexity of  the  Boltzmann entropy ${\rm Ent}$   along Wasserstein geodesics (i.e., $(\ref{Ent(K)})$) is used to characterize the ${\rm CD}(K, \infty)$-condition on Riemannian manifolds $(M, g, \mu)$ 
and metric measure spaces $(X, d, \mu)$, and the inequality $(\ref{S(K, N)})$ for the R\'enyi entropy along geodesics on the Wasserstein space over metric measure space $(X, d, \mu)$ is used to introduce the definition of  the ${\rm CD}(K, N)$-condition on $(X, d, \mu)$. When restricting on Riemannian manifolds with weighted volume measures (i.e., in the setting of smooth metric measure spaces),  we will make a comparison between Lott-Villani and Sturm's definitions of their curvature-dimension ${\rm CD}(K, N)$ conditions ($K\in \mathbb{R}$ and $N\in [n, \infty]$) with our results on the concavity of the Shannon and R\'enyi entropy powers along the geodesic on the Wasserstein space over Riemannian manifolds with the Bakry-Emery ${\rm CD}(K, m)$-condition ($K\in \mathbb{R}$ and $m\in [n, \infty]$) for $N=m$.

In summary, the novelties of our paper are:

(1)  We provide more simple characterizations of ${\rm CD}(K, m)$-condition which are equivalent to the ones of Sturm \cite{Sturm06b} 
and Lott-Villani \cite{LV} on Riemannian manifolds. Our proof uses explicit calculations and is different from \cite{LV, Sturm06b, V2}. Our 
work provides {\it an information-theoretic approach} to better understand the synthetic geometry developed by Lott-Villani and Sturm \cite{LV, V2, Sturm06b}. 

(2) We provide a characterization of the Einstein and quasi-Einstein manifolds by the enhanced entropy differential equality and enhanced entropy power differential equality.  See Theorem \ref{rigidity3}, Theorem \ref{rigidity5}, Theorem \ref{rigidity1} and Theorem \ref{rigidity2} below. These are new in the literature.

 (3) We introduce the Perelman type $W$-entropy associated with the Shannon entropy and the R\'enyi 
 entropy along the Benamou-Brenier geodesic on the Wasserstein space over Riemannian manifolds 
 and  prove its monotonicity and rigidity theorem on the Wasserstein space over complete Riemannian 
 manifolds with CD$(0, m)$-condition. See  Theorem \ref{W-entropy}  and 
 Theorem \ref{W-Rigidity}, which are also new in the literature when $p>1$.  When $p=1$, this was first 
 proved by due to S. Li and the author of this paper in \cite{LL18, LL24}. 

Our work provides a more simple  approach to introduce the CD$(K, N)$-condition on metric measure 
spaces and provides a possible way to define the notion of $n$-dimensional metric measure spaces with 
constant $N$-Ricci curvature, i.e., $Ric_{n, N}=K$, which will be developed in a forthcoming paper in the future.

The rest of this paper is organized as follows. In Section 3, we prove the differential inequalities and rigidity theorem for the Shannon entropy and its power along geodesics on the Wasserstein space over Riemannian manifolds.
  In Section 4, we prove the differential inequalities and rigidity theorem  for the R\'enyi entropy and the 
  associated $W$-entropy  along the Benamou-Brenier geodesics 
on the Wasserstein space over Riemannian manifolds. In Section 5,  we compare our work with the synthetic geometry  developed by Lott-Villani and Sturm and raise some problems for further study in the future.  

\section{Differential inequalities for Shannon entropy}

In this section we prove  the differential inequalities for the Shannon entropy and its power along the geodesic flow on the Wasserstein space over complete Riemannian manifolds with the 
${\rm CD}(K, m)$ curvature-dimension condition.  We prove the NIW formula which indicates the relationship between the Shannon entropy power $N$, the Fisher information $I$ 
and the $W$-entropy associated with the geodesic flow on the Wasserstein space over compact Riemannian manifolds.  
Moreover, we prove that the rigidity models of the enhanced entropy differential inequalities 
are the so-called $(K, m)$-Einstein manifolds, i.e., 
$Ric_{m, n}(L)=Kg$. In particular, when $m=n$, the rigidity models of the enhanced entropy differential inequalities 
are the Einstein manifolds, i.e.,  $Ric=Kg$. 

\subsection{Differential inequalities for Shannon entropy}
Recall the following result which was proved in \cite{LL24} on complete Riemannian manifolds. In compact case, it is well-known, see e.g. \cite{Lo2}.

\begin{theorem}\label{entropygeo}\label{AAA} (\cite{Lo2, LL24}) Let $M$ be a complete Riemannian manifold, and $V\in C^2(M)$. Suppose that $Ric(L)=Ric+\nabla^2V$ is uniformly bounded on $M$, i.e., there exists a constant $C>0$ such that $|Ric(L)|\leq C$. 
Let $(\rho, \phi)$ be smooth solutions to the Benamou-Brenier geodesic which satisfies the following growth conditions
\begin{eqnarray*}
 \int_M \left[|\nabla\log \rho|^2+|\nabla\phi|^2+|\nabla^2\phi|^2+|L\phi|^2+|\nabla L\phi|^2\right]\rho d\mu<\infty,
\end{eqnarray*}
and for a fixed point $o\in M$, some functions $C_i(t)\geq 0$ and $\alpha_i(t)\geq 0$ on $[0, T]$,
\begin{eqnarray*}
C_1(t)e^{-\alpha_1(t) d^2(x, o)}\leq \rho(x, t)\leq C_2(t)e^{\alpha_2 (t)d^2(x, o)}, \ \ \forall x\in M, t\in [0, T],
\end{eqnarray*} 
and
\begin{eqnarray*}
\int_M d^4(x, o)\rho(x, t)d\mu<\infty, \  \ \ \forall t\in [0, T].
\end{eqnarray*}
Then the entropy dissipation formulae hold 
\begin{eqnarray}
{d\over dt } {\rm Ent}(\rho(t))&=&\int_M \nabla \phi \cdot \nabla \rho d\mu=-\int_M L\phi \rho d\mu, \label{ent1}\\
{d^2\over dt^2 }  {\rm Ent}(\rho(t))&=&\int_M \Gamma_2(\phi, \phi) \rho d\mu, \label{ent2}
\end{eqnarray}
where 
\begin{eqnarray*}
\Gamma_2(\phi, \phi):={1\over 2}L|\nabla \phi|^2-\langle \nabla \phi, \nabla L \phi\rangle.
\end{eqnarray*}
Note that, by the generalized Bochner formula  \cite{BE}, we have
\begin{eqnarray*}
\Gamma_2(\phi, \phi)
=\|\nabla^2 \phi\|_{\rm HS}^2+Ric(L)(\nabla \phi, \nabla \phi).
\end{eqnarray*}

\end{theorem}

\begin{theorem}\label{thm2} Under the assumption of Theorem \ref{entropygeo}, suppose that $Ric_{m, n}(L)\geq Kg$, where $m\in [n, \infty)$ and $K\in C(M, \mathbb{R})$.  
Then the following Shannon entropy differential inequality holds
\begin{eqnarray}
-H''\geq {1\over m}H'^{2}+\int_M K |\nabla \phi |^2\rho d\mu+{1\over m}\int_M |L\phi-I|^2\rho d\mu, \label{HHH2}
\end{eqnarray}
and equivalently the Shannon entropy power differential inequality holds
\begin{eqnarray}
{d^2N_m\over dt^2}\leq -{N_m\over m}\int_M K |\nabla \phi|^2\rho d\mu-{N_m\over m^2}\int_M |L \phi-I|^2\rho d\mu,
 \label{N2m5a}
\end{eqnarray}
where 
\begin{eqnarray*}
I=H'(\rho(t))=\int_M L \phi \rho d\mu
\end{eqnarray*}
 is the Fisher information associated with the geodesic flow on $TP_2(M, \mu)$. 
As a consequence, if $Ric_{m, n}(L)\geq Kg$, where $K\in C(M, \mathbb{R})$, then the Shannon entropy differential inequality holds
\begin{eqnarray}
 H''+{1\over m}H'^{2}+\int_M K |\nabla \phi |^2\rho d\mu\leq 0.
\label{HHH2K}
\end{eqnarray}
In particular, when $K$ is a constant, we obtain $(\ref{EDIW})$ and $(\ref{EPIW})$. 
\end{theorem}
{\it Proof}. Note that,
for any $m\geq n$, we have
\begin{eqnarray*}
\Gamma_2( \phi,  \phi)\geq {|L \phi|^2\over m}+Ric_{m, n}(L)(\nabla \phi, \nabla \phi).
\end{eqnarray*}
Thus
\begin{eqnarray*}
-\partial_t^2 H(\rho(t))\geq \int_M\left[{|L \phi|^2\over m}+ Ric_{m, n}(L)(\nabla \phi, \nabla \phi)\right]\rho d\mu.
\end{eqnarray*}
Since $\rho d\mu$ is a probability measure, the Cauchy-Schwarz inequality  yields
\begin{eqnarray*}
\int_M |L \phi|^2\rho d\mu\geq \left(\int_M L \phi \rho d\mu\right)^2=(\partial_t H(\rho(t)))^2,
\end{eqnarray*}
and hence

\begin{eqnarray*}
-\partial_t^2 H(\rho(t))\geq {(\partial_t H(\rho(t)))^2\over m}+ \int_M Ric_{m, n}(L)(\nabla \phi, \nabla \phi)\rho d\mu.
\end{eqnarray*}
Thus, if $Ric_{m, n}(L)\geq K$ for some function $K\in C(M, \mathbb{R})$, we have
\begin{eqnarray*}
-\partial_t^2 H(\rho(t))\geq {(\partial_t H(\rho(t)))^2\over m}+  \int_M K |\nabla \phi|^2\rho d\mu.
\end{eqnarray*}
This finishes the proof of $(\ref{HHH2K})$. In particular, when $K$ is a constant, we obtain $(\ref{EDIW})$. 

Note that, along the geodesic flow $(\rho(t), \phi(t))$ on $T\mathcal{P}_2(M, \mu)$, we have
\begin{eqnarray*}
\int_M |\nabla \phi(t)|^2\rho(t)d\mu=W_2^2(\rho(0), \rho(1)), \ \ \forall t>0.  \label{W2}
\end{eqnarray*}
This proves the equivalent form of the entropy differential inequality. 

On the other hand,  the Shannon entropy power satisfies
\begin{eqnarray*}
{d\over dt}N_m(\rho(t))&=&N_m(\rho(t)){\partial_t H(\rho(t))\over m}, \label{NH1}\\
{d^2\over dt^2}N_m(\rho(t))&=&N_m(\rho(t))\left[{\partial^2_t H(\rho(t))\over m}+\left({\partial_t H(\rho(t))\over m}\right)^2\right]. \label{NH2}
\end{eqnarray*}
Thus the EPDI $(\ref{EPIW})$  follows from EDI $(\ref{EDIW})$.  \hfill $\square$

Let
\begin{eqnarray*}
\Theta_K(\rho(t))=\int_M K|\nabla \phi(t)|^2\rho(t) d\mu.
\end{eqnarray*}
Then we have the Riccati differential inequality along $(\rho(t), \phi(t))$
\begin{eqnarray}
H''+{1\over m}H'^{2}+\Theta_K(\rho(t))\leq 0. 
\label{Riccati-1}
\end{eqnarray}
The boundary condition for $H'(\rho(t))$ is given by
\begin{eqnarray*}
H'(0+)=\lim\limits_{t\rightarrow 0+}H'(\rho(t)).
\end{eqnarray*}Note that $(\ref{Riccati-1})$ is equivalent to 
\begin{eqnarray*}
{d^2N_m\over dt^2}\leq -{\Theta_K(\rho(t)) \over m}N_m.
 \label{N2m5Kaa}
\end{eqnarray*}
Let $H_{m, K}'(t)$ be the unique solution to the Riccati equation 
\begin{eqnarray}
H_{m, K}''+{1\over m}H_{m, K}'^2+\Theta_K(\rho(t))=0
\label{Riccati-2}
\end{eqnarray}
with the initial value $H_{m, K}'(0)\geq H'(0+)$, and let $N_{m, K}$ be the unique solution to the second order differential equation 
\begin{eqnarray}
{d^2N_{m, K}\over dt^2}=-{\Theta_K(\rho(t))\over m}N_{m, K}
 \label{N2m5Kma}
\end{eqnarray}
with two-point boundary condition
$$N_{m, K}(0)=N_m(\rho(0)), \ \ \ N_{m, K}(1)=N_{m}(\rho(1)).$$

By the Sturm-Liouville comparison theorem and the maximum principle for Riccati equation, we have

\begin{theorem}\label{thm2b} Under the same condition as in Theorem \ref{thm2}, we have
\begin{eqnarray}
H'(\rho(t))\leq H_{m, K}'(t), \label{Riccati-3}
\end{eqnarray}
and
\begin{eqnarray}
N_{m}(\rho(t))\geq N_{m, K}(t).
 \label{N2m5KmaK}
\end{eqnarray}
In the case $K\in \mathbb{R}$ is a constant, we have $\Theta_K(\rho(t))=K\theta:=KW_2^2(\rho(0)\mu, \rho(1)\mu)$, and
$(\ref{Riccati-1})$ and $(\ref{Riccati-2})$ read as follows
\begin{eqnarray}
H_{m, K}''+{1\over m}H_{m, K}'^2+K\theta^2=0, \label{Riccati-2K}
\end{eqnarray}
and 
\begin{eqnarray}
{d^2N_{m, K}\over dt^2}=-{K\theta^2\over m}N_{m, K}. \label{N2m5KmaK}
\end{eqnarray}
Hence $N_{m, L}(t)$ is given by $(\ref{NmKk})$, i.e., 
\begin{eqnarray*}
N_{m, K}(t)=\sigma_{K, m, \theta}(1-t)N_m(\rho(0))+ \sigma_{K, m, \theta}(t)N_m(\rho(1)). 
\end{eqnarray*}
\end{theorem}

{\it Proof}. 
Note that $(\ref{Riccati-1})$ and $(\ref{Riccati-2})$ have the constant term $\Theta_K(\rho(t))=\int_M K|\nabla \phi(t)|^2\rho(t) d\mu$ which depends on $(\rho(t), \phi(t))$. Let $u=H'(t)-H_{m, K}'(t)$. We have
\begin{eqnarray*}
u'\leq -{H'+H_{m, K}'\over m}u.
\label{Riccati-5}
\end{eqnarray*}
By the Gronwall inequality, we have
\begin{eqnarray*}
u(t)\leq u(0)\exp\left(-{H(t)+H_{m, K}(t)-H(0)-H_{m, K}(0)\over m}\right).
\label{Riccati-6}
\end{eqnarray*}
In particular, as $u(0)=H'(t)-H_{m, K}'(t)\leq 0$, we have
\begin{eqnarray*}
u(t)\leq 0.
\label{Riccati-7}
\end{eqnarray*}
Integrating along geodesic from $\rho(0)$ to $\rho(1)$, the Shannon entropy power inequality $(\ref{EPIW})$ implies that
\begin{eqnarray}
N_{m}(\rho(t))\geq N_{m, K}(t).
 \label{N2m5KmaK}
\end{eqnarray}
From which we derive  $(\ref{NmKk})$. 
\hfill $\square$

%

\subsection{The $W$-entropy and NIW formula for Shannon entropy}

Recall that, in \cite{LL18, LL24}, S. Li and the author of this paper introduced the $W$-entropy and proved the following $W$-entropy formula for the Shannon entropy along
 the geodesic flow on the Wasserstein space $P_2(M, \mu)$.

\begin{theorem}\label{MT2} \cite{LL18, LL24}
Let $(M, g)$  be a complete Riemannian manifold, $V\in C^2(M)$ such that $Ric(L)=Ric+\nabla^2 V$ 
is uniformly bounded. 
Let $(\rho, \phi): M\times [0, T]\rightarrow \mathbb{R}^+\times \mathbb{R}$  be a smooth solution to 
the Benamou-Brenier geodesic 
on $P_2(M, \mu)$ satisfying the growth condition as required in Theorem \ref{AAA}. For any $m\geq n$, define the $H_m$-entropy and $W_m$-entropy for the geodesic flow $(\rho, \phi)$ on $TP^\infty_2(M, \mu)$ as follows
\begin{eqnarray*}
H_m(\rho, t)=H(\rho(t))-{m\over 2}\left(1+\log(4\pi t^2)\right),
\end{eqnarray*}
and
\begin{eqnarray*}
W_m(\rho, t)={d\over dt}(tH_m(\rho, t)).
\end{eqnarray*}
Then for all $t>0$, we have
\begin{eqnarray}
{1\over t}{d\over dt}W_m(\rho, t)&=&-\int_M \left[\left\|\nabla^2 \phi-{g\over t}\right\|_{\rm HS}^2+Ric_{m,
n}(L)(\nabla \phi, \nabla \phi) \right]\rho d\mu\nonumber\\
& &\ \ \ \ \ \ \ \ \ \ -{1 \over m-n}\int_M \left|\nabla V\cdot
\nabla \phi-{m-n\over t}\right|^2 \rho d\mu.\label{Wgeo}
\end{eqnarray}
In paricular, if $Ric_{m, n}(L)\geq 0$, then $W_m(\rho, t)$ is nonincreasing in time $t$ along the geodesic flow on $TP^\infty_2(M, \mu)$.  
\end{theorem}

Moreover,   S. Li and the author \cite{LL18, LL24} proved  the following rigidity theorem of the  $W$-entropy for the Shannon entropy along the Benamou-Brenier geodesic on the Wasserstein space over Riemannian manifold with CD$(0, m)$-condition. 

\begin{theorem} \label{W-Rigidity Shannon} \cite{LL18, LL24} Let $(M, g, V)$ be a complete 
Riemannian manifold with $Ric(L)=Ric+\nabla^2V$ uniformly bounded, and let $(\rho(t), \phi(t))$ be a smooth Benamou-Brenier geodesic on $P_2(M, \mu)$ satisfying  the growth condition as required in Theorem \ref{AAA}. Suppose that 
$Ric_{m, n}(L)\geq 0$ and ${d\over dt}W_{m}(\rho, t)=0$ at some $t=\tau>0$. Then $M$ is isometric to $\mathbb{R}^n$, $m=n$, $V$ is a constant and $(\rho, \phi)=(\rho_n, \phi_n)$. 
\end{theorem}

As a consequence of the $W$-entropy formula, we can derive the following $W$-entropy inequality: On complete Riemannian manifold with bounded geometry condition, if $Ric_{m, n}(L)\geq K$, then we have
\begin{eqnarray*}
{1\over t}{d\over dt}W_m(\rho, t)\leq-\int_M K|\nabla \phi|^2\rho d\mu-{1\over m}\int_M \left|\Delta\phi-{m\over t}\right|^2\rho d\mu. \label{Wgeo-2}
\end{eqnarray*}
In particular, for $L=\Delta$, $m=n$ and $\mu=v$, we obtain: on any complete Riemannian manifold with bounded geometry condition, if $Ric\geq K$, we have
\begin{eqnarray*}
{1\over t}{d\over dt}W_n(\rho, t)\leq-\int_M K|\nabla \phi|^2\rho dv-{1\over n}\int_M \left|\Delta\phi-{n\over t}\right|^2\rho dv. \label{Wgeo-3}
\end{eqnarray*}

In this section, we prove the following NIW formula, which indicates an interesting relationship between the Shannon entropy power $N$, the Fisher information $I$ and the $W$-entropy associated with the geodesic flow on the Wasserstein space over (weighted) complete Riemannian manifolds.  More precisely, we have the following

%
%

\begin{theorem}\label{thm5} Let $M$ be an $n$-dimensional complete Riemannian manifold with $V\in C^2(M)$. 
Define
\begin{eqnarray*}
H(\rho(t))=-\int_M \rho \log \rho d\mu, \ \ N_m(\rho(t))=e^{H(\rho(t))\over m}.
\end{eqnarray*}
Then 
\begin{eqnarray}
{m\over N_m} {d^2N_m\over dt^2}&=&-{1\over m} \int_M \left| L \phi- \int_M L \phi  \rho d\mu \right|^2 \rho d\mu-
\int_M Ric_{m, n}(L)(\nabla \phi, \nabla \phi) \rho d\mu\nonumber\\
& &-{m-n\over mn} \int_M \left(L \phi+{m\over m-n}\nabla V\cdot\nabla \phi\right)^2\rho d\mu- \int_M \left\|\nabla^2 \phi-{\Delta \phi\over n} g\right\|^2_{\rm HS}\rho d\mu. \nonumber \\
& &\label{N2m1}
\end{eqnarray}
Moreover, the following NIW formula holds
\begin{eqnarray*}
 {d^2N_m\over dt^2}={ N_m\over m} \left[{1\over m}\left|I-{m\over t}\right|^2+{1\over t}{dW_m\over dt}\right]. \label{N2m1b}
\end{eqnarray*}
where $I=H'(\rho(t))=\int_M L\phi(t) \rho(t) d\mu$ is the Fisher information. 
In particular, if $Ric_{m, n}(L)\geq K$, we have
\begin{eqnarray}
 {1\over t}{dW_m\over dt}\leq -\int_M K|\nabla \phi|^2\rho d\mu-{1\over m}\left|I-{m\over t}\right|^2. \label{IKm-1} 
 \end{eqnarray}
\end{theorem}
{\it Proof}. Using $\|A\|_{\rm HS}^2={|{\rm Tr A}|^2\over n}+\left\|A-{{\rm Tr}A\over n} g\right\|^2_{\rm HS}$, we have

\begin{eqnarray*}
\|\nabla^2 \phi\|_{\rm HS}^2= {|\Delta \phi|^2\over n}+\left\|\nabla^2 \phi-{\Delta \phi\over n} g\right\|^2_{\rm HS}.
\end{eqnarray*}
Applying the elementary equality
\begin{eqnarray*}
(a+b)^2={a^2\over 1+\varepsilon}-{b^2\over \varepsilon}+{\varepsilon\over 1+\varepsilon}\left(a+{1+\varepsilon\over \varepsilon}b\right)^2
\end{eqnarray*}
to $a=L \phi$, $b=\nabla V\cdot\nabla \phi$ and $\varepsilon={m-n\over n}$, we have
\begin{eqnarray*}
|\Delta \phi|^2= {n\over m}|L \phi|^2-{n\over m-n}|\nabla V\cdot\nabla \phi|^2+{m-n\over m}\left(L \phi+{m\over m-n}\nabla V\cdot\nabla \phi\right)^2,
\end{eqnarray*}
and 
\begin{eqnarray*}
\|\nabla^2 \phi\|_{\rm HS}^2&=&  {|L \phi|^2\over m}-{|\nabla V\cdot\nabla \phi|^2\over m-n}+{m-n\over mn}\left(L \phi+{m\over m-n}\nabla V\cdot\nabla \phi\right)^2\\
& &\hskip2cm +\left\|\nabla^2 \phi-{\Delta \phi\over n} g\right\|^2_{\rm HS}.
\end{eqnarray*}
This yields
\begin{eqnarray*}
\Gamma_2( \phi,  \phi)
&=& {|L \phi|^2\over m}+Ric_{m, n}(L)({\nabla \phi, \nabla \phi)}\\
& &+{m-n\over mn}\left(L \phi+{m\over m-n}\nabla V\cdot\nabla \phi\right)^2 +\left\|\nabla^2 \phi-{\Delta \phi\over n} g\right\|^2_{\rm HS}.
\end{eqnarray*}
Substituting this into $(\ref{ent2})$, we have
\begin{eqnarray*}
\partial_t H(\rho(t))&=&\int_M L \phi\rho d\mu,\\
\partial_t ^2 H(\rho(t))&=&-\int_M \left[{|L \phi|^2\over m}+Ric_{m, n}(L)({\nabla \phi, \nabla \phi)}\right]\rho d\mu\\
& &-{m-n\over mn}\int_M \left(L \phi+{m\over m-n}\nabla V\cdot\nabla \phi\right)^2\rho d\mu-\int_M \left\|\nabla^2 \phi-{\Delta \phi\over n} g\right\|^2_{\rm HS}\rho d\mu.
\end{eqnarray*}
Combining this with $(\ref{ent1})$, we have

\begin{eqnarray}
& &\partial^2_t H(\rho(t))+{(\partial_t H(\rho(t))^2\over m}\nonumber\\
&=&-{1\over m} \int_M \left| L \phi- \int_M L \phi  \rho d\mu \right|^2 \rho d\mu-\int_M Ric_{m, n}(L)(\nabla \phi, \nabla \phi) \rho d\mu\nonumber\\
& &-{m-n\over mn} \int_M \left(L \phi+{m\over m-n}\nabla V\cdot\nabla \phi\right)^2\rho d\mu- \int_M \left\|\nabla^2 \phi-{\Delta \phi\over n} g\right\|^2_{\rm HS}\rho d\mu. \label{HHHmK}
\end{eqnarray}
Combining this with $(\ref{NH1})$ and $(\ref{NH2})$, we have

\begin{eqnarray}
{m\over N_m} {d^2N_m\over dt^2}
&=&-{1\over m} \int_M \left| L \phi- \int_M L \phi  \rho d\mu \right|^2 \rho d\mu-\int_M Ric_{m, n}(L)(\nabla \phi, \nabla \phi) \rho d\mu\nonumber\\
& &-{m-n\over mn} \int_M \left(L \phi+{m\over m-n}\nabla V\cdot\nabla \phi\right)^2\rho d\mu- \int_M \left\|\nabla^2 \phi-{\Delta \phi\over n} g\right\|^2_{\rm HS}\rho d\mu.\nonumber\\
\label{qqq}
\end{eqnarray}
By  $(\ref{qqq})$ and $(\ref{Wgeo})$, we have
\begin{eqnarray*}
& &{m\over N_m} {d^2N_m\over dt^2}-{1\over t}{dW_m\over dt}\\
&=&-{1\over m} \int_M \left| L \phi- \int_M L \phi  \rho d\mu \right|^2 \rho d\mu-{m-n\over mn} \int_M \left(L \phi+{m\over m-n}\nabla V\cdot\nabla \phi\right)^2\rho d\mu\\
& &- \int_M \left\|\nabla^2 \phi-{\Delta \phi\over n} g\right\|^2_{\rm HS}\rho d\mu+\int_M \left\|\nabla^2 \phi-{g\over t}\right\|_{\rm HS}^2 \rho d\mu +{1\over m-n}
\int_M \left|\nabla V\cdot\nabla \phi+{m-n\over t}\right|^2\rho d\mu.
\end{eqnarray*}
Note that
\begin{eqnarray*}
\left\|\nabla^2 \phi-{g\over t}\right\|_{\rm HS}^2
&=&{1\over n} \left|\Delta \phi-{n\over t}\right|^2+\left\|\nabla^2 \phi-{\Delta \phi\over n}g\right\|_{\rm HS}^2\\
&=&{1\over n} \left|L \phi-{m\over t}+\left(\nabla V\cdot\nabla \phi+{m-n\over t}\right)   \right|^2+\left\|\nabla^2 \phi-{\Delta \phi\over n}g\right\|_{\rm HS}^2\\\
&=& {1\over m} \left|L \phi-{m\over t}  \right|^2-{1\over m-n}\left(\nabla V\cdot\nabla \phi+{m-n\over t}\right)^2\\
& &+{m-n\over mn}\left[L \phi+{m\over m-n}\nabla V\cdot\nabla \phi\right]^2+\left\|\nabla^2 \phi-{\Delta \phi\over n}g\right\|_{\rm HS}^2.
\end{eqnarray*}
Combining this with the fact $I=\partial_t  H(\rho(t))=\int_M L \phi \rho d\mu$, we have
\begin{eqnarray*}
{m\over N_m} {d^2N_m\over dt^2}={1\over t}{dW_m\over dt}+{1\over m}\left|I-{m\over t}\right|^2.
\end{eqnarray*}
This proves the NIW formula. The $W$-inequality $(\ref{IKm-1})$ follows from $(\ref{N2m5Kaa})$. \hfill $\square$

\medskip

In particular, when $m=n$, $V=0$ and $L=\Delta$, we have the following

\begin{theorem}\label{thm6} Let $M$ be a complete Riemannian manifold. Then 
\begin{eqnarray}
{n\over N_n} {d^2N_n\over dt^2}&=&-{1\over n} \int_M \left| \Delta \phi- \int_M \Delta \phi  \rho dv \right|^2 \rho dv
-\int_M Ric(\nabla \phi, \nabla \phi) \rho dv\nonumber\\
& &\hskip2cm - \int_M \left\|\nabla^2 \phi-{\Delta \phi\over n} g\right\|^2_{\rm HS}\rho dv.\label{N2n1}
\end{eqnarray}
Moreover, the following NIW formula holds
\begin{eqnarray}
 {d^2N_n\over dt^2}={ N_n\over n} \left[{1\over n}\left|I-{n\over t}\right|^2+{1\over t}{dW_n\over dt}\right].\label{N2n2}
\end{eqnarray}
In particular, if $Ric\geq K$, we have
\begin{eqnarray}
 {1\over t}{dW_n\over dt}\leq -\int_M K|\nabla \phi|^2\rho dv-{1\over n}\left|I-{n\over t}\right|^2. \label{IKn-2} 
 \end{eqnarray}

\end{theorem}

\begin{remark} Using again $\|A\|_{\rm HS}^2={|{\rm Tr A}|^2\over n}+\left\|A-{{\rm Tr}A\over n} g\right\|^2_{\rm HS}$, we have
\begin{eqnarray*}
\left\|\nabla^2 \phi-{I\over n}g\right\|_{\rm HS}^2&=&{1\over n}\left|\Delta \phi-I\right|^2+\left\|\nabla^2 \phi-{I\over n}g-{1\over n}\left(\Delta \phi-I\right)g\right\|_{\rm HS}^2\nonumber\\
&=&{1\over n}\left|\Delta \phi-I\right|^2+\left\|\nabla^2 \phi-{\Delta \phi\over n}g\right\|_{\rm HS}^2.
\end{eqnarray*}
This yields another explicit expression of the second derivative of the Shannon entropy power along the geodesic flow on $TP_2(M, v)$, 
which is indeed equivalent to $(\ref{N2n1})$. More precisely, we have
\begin{eqnarray*}
{d^2 N_n\over dt^2}
&=&-{N_n\over n}\int_M \left(Ric(\nabla \phi, \nabla \phi)+\left\|\nabla^2 \phi-{I\over n}g\right\|_{\rm HS}^2\right) \rho dv.\label{N2n3}
\end{eqnarray*}
Similarly, for $m>n$, $I=\int_M L\phi \rho d\mu$, using
\begin{eqnarray*}
\left\|\nabla^2 \phi-{I\over m}g\right\|_{\rm HS}^2&=&{1\over n}\left|\Delta \phi-{n\over m}I\right|^2+\left\|\nabla^2 \phi-{I\over m}g-{1\over n}\left(\Delta \phi-{n\over m}I\right)g\right\|_{\rm HS}^2\\
&=&{1\over n}\left|\Delta \phi-{n\over m}I\right|^2+\left\|\nabla^2 \phi-{\Delta \phi\over n}g\right\|_{\rm HS}^2,
\end{eqnarray*}
and
\begin{eqnarray*}
{1\over n}|\Delta \phi-{n\over m}I|^2&=& {1\over m}|L \phi-I|^2-{1\over m-n}\left|\nabla V\cdot\nabla \phi+{m-n\over m}\right|^2\\
& &+{m-n\over mn}\left(L \phi+{m\over m-n}\nabla V\cdot\nabla \phi\right)^2,
\end{eqnarray*}
we have

\begin{eqnarray*}
 {m\over N_m}{d^2N_m\over dt^2}&=&-\int_M \left(Ric_{m, n}(L)(\nabla \phi, \nabla \phi) + \left\|\nabla^2 \phi-{I\over m} g\right\|^2_{\rm HS}\right)\rho d\mu\\
 & &-{1\over (m-n)} \int_M \left(\nabla V\cdot\nabla \phi+{m-n\over m}\right)^2\rho d\mu.
\end{eqnarray*}

\end{remark}

\subsection{Rigidity of Shannon entropy and entropy power}

In this subsection, we prove the rigidity theorems for the  Shannon entropy and entropy power inequality 
under the ${\rm CD}(K, m)$-condition. 

\begin{theorem} \label{rigidity1} Under the same condition and notation as  in Theorem \ref{Thm1}, we have the 
enhanced entropy power differential inequality
\begin{eqnarray}
{d^2N_n\over dt^2}+{N_n\over n^2}\int_M |\Delta \phi-I|^2\rho dv\leq -{KN_n\over n}\int_M |\nabla \phi|^2\rho dv.
 \label{N2n5a}
\end{eqnarray}
This indeed improves the entropy power concavity inequality in Theorem \ref{Thm1}, i.e., 
\begin{eqnarray}
{d^2N_n\over dt^2}\leq -{KN_n\over n}\int_M |\nabla \phi|^2\rho dv.
 \label{N2n6}
\end{eqnarray}
{Moreover, the equality in the enhanced entropy power differential inequality $(\ref{N2n5a})$  holds for all geodesics $(\rho(t), \phi(t))$ if and only if $(M, g)$ is an Einstein manifold with constant Ricci 
curvature $Ric=Kg$. 
}
\end{theorem}
{\it Proof}.  By $(\ref{N2n1})$ in Theorem \ref{thm6} or $(\ref{N2n3})$,  under the condition $Ric\geq Kg$, we have 
\begin{eqnarray}
{d^2N_n\over dt^2}&\leq& -{KN_n\over n}\int_M |\nabla \phi|^2\rho dv-{N_n\over n}\int_M \left\|\nabla^2 \phi-{I\over n} g\right\|^2_{\rm HS}\rho dv.\label{N2n4}
\end{eqnarray}
Using the trace inequality $\left\|\nabla^2 \phi-{I\over n} g\right\|^2_{\rm HS}\geq {1\over n}|\Delta \phi-I|^2$, we have
\begin{eqnarray}
{d^2N_n\over dt^2}\leq -{KN_n\over n}\int_M |\nabla \phi|^2\rho dv-{N_n\over n^2}\int_M |\Delta \phi-I|^2\rho dv.
 \label{N2n5}
\end{eqnarray}
This improves the entropy power concavity inequality in Theorem \ref{Thm1}, i.e., $(\ref{EPIW})$. 

{
Moreover, $(\ref{N2n5a})$ is an equality if and only if 
\begin{eqnarray}
Ric(\nabla \phi, \nabla \phi)(x)=K|\nabla \phi(x)|^2.
\label{N2n7}
\end{eqnarray}
}
{Since the set $\{\nabla \phi(x): \phi\in C^\infty(M)\}$ spans the whole tangent space $T_xM$, we conclude that $(M, g)$ is Einstein manifold with constant Ricci 
curvature $Ric=Kg$. 
 }\hfill $\square$
\medskip

\begin{remark} In particular, if $Ric\geq 0$, then $I(t)=H'(t)$ satisfies
$I'+{I^2\over n}\leq 0$. This yields ${d\over dt}\left({1\over I}-{t\over n}\right)\geq 0$, and hence
$I(t)\leq {n\over Cn+t}$, 
where 
$C:=\lim\limits_{t\rightarrow 0+}\left({1\over I}-{t\over n}\right)={1\over I(0)}.$
\end{remark}

In particular, we have the following 

\begin{theorem}  
Under the assumption $Ric\geq 0$ and $p\geq 1-{1\over n}$,  we have the 
enhanced entropy power differential inequality
\begin{eqnarray}
{d^2N_n\over dt^2}+{N_n\over n^2}\int_M |\Delta \phi-I|^2\rho dv\leq 0.
 \label{N2n5a=0}
\end{eqnarray}In particular, we have
\begin{eqnarray}
{d^2\over dt^2} N_{n}(\rho(t))\leq 0, \label{Nn0}
\end{eqnarray}
i.e.,  $N_{n}(\rho(t))$ is concave in $t$ on $[0, \infty)$. 
Moreover, the equality in the enhanced entropy power differential inequality $(\ref{N2n5a=0})$  holds for all geodesics $(\rho(t), \phi(t))$ if and only if $(M, g)$ is Ricci flat, i.e., $Ric=0$.  

Furthermore, in the case where $Ric=0$, the 
equality in \eqref{Nn0} holds if and only if $\phi$ satisfies the Hessian soliton equation
\begin{eqnarray*}
\nabla^2 \phi={I\over n}g.
\end{eqnarray*}
In the case $(M, g)$ is a complete Riemannian manifold with bounded geometry condition as in Theorem \ref{AAA}, and $(\rho(t), \phi(t))$ is a smooth Benamou-Brenier geodesic on the Wasserstein space over Riemannian manifold satisfying  growth condition as required in Theorem \ref{AAA}, then under the 
assumption $Ric=0$ and $p\geq 1-{1\over n}$,  the equality  in $(\ref{Nn0})$ holds at some $t=\tau>0$ if and only if $M$ is isometric to $\mathbb{R}^n$, and $(\rho, \phi)=(\rho_n, \phi_n)$.
\end{theorem}
{\it Proof}. The last statement of rigidity theorem is indeed a reformulation of Theorem \ref{W-Rigidity Shannon}. \hfill $\square$

\medskip
In general case of weighted Riemannian manifolds, we have the following

\begin{theorem} \label{rigidity2} Under the same condition and notation as  in Theorem \ref{thm2}, we have the 
enhanced entropy power differential inequality
\begin{eqnarray}
{d^2N_m\over dt^2}d\mu+{N_m\over m^2}\int_M |L \phi-I|^2\rho d\mu\leq -{KN_m\over m}\int_M |\nabla \phi|^2\rho d\mu.
 \label{N2m5}
\end{eqnarray}
This indeed improves the entropy power concavity inequality EPCI in Theorem \ref{thm2}, i.e., 
\begin{eqnarray}
{d^2N_m\over dt^2}\leq -{KN_m\over m}\int_M |\nabla \phi|^2\rho d\mu.
 \label{N2m6}
\end{eqnarray}
{Moreover, the equality in the enhanced entropy power differential inequality $(\ref{N2m5})$  holds for all geodesics $(\rho(t), \phi(t))$ if and only if $(M, g, V)$ is a $(K, m)$-Einstein manifold with constant $m$-dimensional Bakry-Emery Ricci 
curvature  
\begin{eqnarray*}
Ric_{m, n}(L)=Kg.
\end{eqnarray*}
}
\end{theorem}
{\it Proof}.  By the NIW formula $(\ref{N2m1})$ in Theorem \ref{thm5}, under the condition $Ric_{m, n}(L)\geq Kg$, where $K\in C(M, \mathbb{R})$, we have 
\begin{eqnarray}
{m\over N_m} {d^2N_m\over dt^2}&\leq &-\int_M K|\nabla \phi|^2 \rho d\mu-{m-n\over mn} \int_M \left(L \phi+{m\over m-n}\nabla V\cdot\nabla \phi\right)^2\rho d\mu\nonumber\\
& &-{1\over m} \int_M \left| L \phi- \int_M L \phi  \rho d\mu \right|^2 \rho d\mu- \int_M \left\|\nabla^2 \phi-{\Delta \phi\over n} g\right\|^2_{\rm HS}\rho d\mu. \label{N2m3}
\end{eqnarray}
This yields 
\begin{eqnarray}
{d^2N_m\over dt^2}\leq -{N_m\over m}\int_M K |\nabla \phi|^2\rho d\mu-{N_m\over m^2}
\int_M |L \phi-I|^2\rho d\mu, \label{N2m4}
\end{eqnarray}
which implies 
\begin{eqnarray}
{d^2N_m\over dt^2} \leq -{N_m\over m}\int_M K |\nabla \phi|^2\rho d\mu. \label{N2m5b}
\end{eqnarray}
In particular, when $K$ is a constant, we have 
\begin{eqnarray}
{m\over N_m} {d^2N_m\over dt^2}&\leq &-K\int_M |\nabla \phi|^2 \rho d\mu-{m-n\over mn} \int_M \left(L \phi+{m\over m-n}\nabla V\cdot\nabla \phi\right)^2\rho d\mu\nonumber\\
& &-{1\over m} \int_M \left| L \phi- \int_M L \phi  \rho d\mu \right|^2 \rho d\mu- \int_M \left\|\nabla^2 \phi-{\Delta \phi\over n} g\right\|^2_{\rm HS}\rho d\mu. \label{N2m3b}
\end{eqnarray}
This yields 
\begin{eqnarray}
{d^2N_m\over dt^2}\leq -{KN_m\over m}\int_M |\nabla \phi|^2\rho d\mu-{N_m\over m^2}\int_M |L \phi-I|^2\rho d\mu.
 \label{N2m5K=c}
\end{eqnarray}

The inequalities $(\ref{N2m4})$, $(\ref{N2m5})$ and $(\ref{N2m5K=c})$ are improved versions of the entropy power concavity inequality in Theorem \ref{thm2}, i.e., if $Ric_{m, n}(L)\geq K$ for a constant 
$K\in \mathbb{R}$, then
\begin{eqnarray}
{d^2N_m\over dt^2}&\leq &-{KN_m\over m}\int_M |\nabla \phi|^2 \rho d\mu. \label{N2m6K=c}
\end{eqnarray}

{Moreover,  we can conclude that, the inequality $(\ref{N2m3b})$ 
becomes an equality for all geodesics $(\rho(t), \phi(t))$  if and only if 
\begin{eqnarray*}
Ric_{m, n}(L)(\nabla \phi, \nabla \phi)=K|\nabla \phi|^2. 
\end{eqnarray*}
Since $\nabla \phi$ spans the whole tangent space $TM$, we conclude that $(M, g)$ is a $(K, m)$-Einstein manifold with constant $m$-dimensional Bakry-Emery Ricci 
curvature 
\begin{eqnarray*}
Ric_{m, n}(L)=Kg.
\end{eqnarray*}
This completes the proof of Theorem \ref{thm2}.} \hfill $\square$

\begin{remark}
{In the case where $Ric_{m, n}(L)=Kg$, we can check that $(\ref{N2m5b})$ becomes an equality if and only if the differential inequality $(\ref{HHH2K})$ becomes an equality, i.e., 
\begin{eqnarray}
-H''(\rho(t))={H'^2(\rho(t))\over m}+\int_M K|\nabla \phi(t)|^2\rho(t) d\mu. \label{EDIW1Kb}
\end{eqnarray}
In particular, the inequality $(\ref{N2m6K=c})$ becomes an equality if and only if the differential inequality $(\ref{EDIW})$ becomes the equality  \eqref{Riccati-2K}, i.e., 
\begin{eqnarray}
H''(\rho(t))+{H'^2(\rho(t))\over m}+KW_2^2(\rho(0), \rho(1))=0. \label{EDIW1bK=c}
\end{eqnarray}
Equivalently, the Fisher information $I=H'(\rho(t))$ satisfies the Riccati equation
\begin{eqnarray}
I'(\rho(t))+{I^2(\rho(t))\over m}+KW_2^2(\rho(0), \rho(1))=0. \label{Riccatti-I}
\end{eqnarray}
}
\end{remark}

\begin{remark}
In particular, if  $Ric_{m, n}(L)\geq 0$, then  $I(t)=H'(t)$ satisfies
$I'+{I^2\over m}\leq 0.$ This yields 
${d\over dt}\left({1\over I}-{t\over m}\right)\geq 0$, and hence $I(t)\leq {m\over t+Cm}$, 
where 
$C:=\lim\limits_{t\rightarrow 0}\left({1\over I}-{t\over m}\right)={1\over I(0)}.$

\end{remark}

In particular, we have the following 

\begin{theorem} Under the assumption $Ric_{m, n}(L)\geq 0$ and $p\geq 1-{1\over m}$,  we have the 
enhanced entropy power differential inequality\begin{eqnarray}
{d^2N_m\over dt^2}d\mu+{N_m\over m^2}\int_M |L \phi-I|^2\rho d\mu\leq 0.
 \label{N2m5K=0}
\end{eqnarray}
In particular, we have
\begin{eqnarray}
{d^2\over dt^2} N_{m}(\rho(t))\leq 0, \label{Nm0}
\end{eqnarray}
i.e.,  $N_{m}(\rho(t))$ is concave in $t$ on $[0, \infty)$.  
Moreover, under the assumption $Ric_{m, n}(L)\geq 0$ and $p\geq 1-{1\over m}$,  the equality  in 
$(\ref{N2m5K=0})$ 
holds for all geodesics $(\rho, \phi)$ at some $t=\tau>0$ if and only if $(M, g, V)$ is quasi-Ricci flat, i.e., $Ric_{m, n}(L)=0$.  Furthermore, in the case where $Ric_{m, n}(L)=0$, the equality  in 
$(\ref{Nm0})$ holds if and only if $\phi$ satisfies the Hessian soliton equation
\begin{eqnarray*}
\nabla^2 \phi={I\over m}g.
\end{eqnarray*}
In the case $(M, g, V)$ is a complete Riemannian manifold with bounded geometry condition as in Theorem \ref{AAA}, and $(\rho(t), \phi(t))$ is a smooth Benamou-Brenier geodesic on the Wasserstein space over Riemannian manifold satisfying  growth condition as required in Theorem \ref{AAA}. Then, under the assumption $Ric_{m, n}(L)=0$ and $p\geq 1-{1\over m}$,  the equality  in $(\ref{Nm0})$ holds at some $t=\tau>0$ if and only if $M$ is isometric to $\mathbb{R}^n$, $m=n$, $V$ is a constant and $(\rho, \phi)=(\rho_n, \phi_n)$.
\end{theorem}
{\it Proof}. The last statement of rigidity theorem is indeed a reformulation of Theorem \ref{W-Rigidity Shannon}. \hfill $\square$

\section{Differential  inequalities for R\'enyi entropy}

In this section we prove  the differential inequalies for the R\'enyi entropy and its power along the geodesic flow on the Wasserstein space over complete Riemannian manifolds with the 
${\rm CD}(K, m)$ curvature-dimension condition.    We prove the NIW formula which indicates the relationship between the R\'enyi entropy power $N$, the Fisher information $I$ 
and the $W$-entropy associated with the geodesic flow on the Wasserstein space over complete Riemannian manifolds.   
Moreover, we prove that the rigidity models of the R\'enyi entropy power $N$ is the $(K, m)$-Einstein manifolds, i.e., 
$Ric_{m, n}(L)=Kg$. In particular, when $m=n$, the rigidity models of the R\'enyi entropy power $N$ is the $K$-Einstein manifolds, i.e.,  $Ric=Kg$. 

Let $(M, g)$ be an $n$-dimensional complete Riemannian manifold with a weighted volume measure $d\mu=e^{-V}dv$ satisfying the Bakry-Emery curvature-dimension ${\rm CD}(K, m)$-condition, i.e.,  $Ric_{m, n}(L)\geq Kg$, where $m\in \mathbb{R}^+$ and $K\in \mathbb{R}$, $m\geq n$. 

Let $(\rho, \phi)$ be a smooth solution to the  geodesic flow on $T{P}_2(M, \mu)$, i.e., $(\rho, \phi)$ is a smooth solution to the continuity equation $(\ref{TA})$ and the Hamilton-Jacobi equation $(\ref{HJ})$. 
Define the $p$-th R\'enyi entropy by 
\begin{eqnarray}
H_p(\rho)={1\over 1-p}\log \int_M \rho^p d\mu,
\end{eqnarray} and the $(m, p)$-R\'enyi entropy power $N_{m, p}$ by 
\begin{eqnarray}
N_{m, p}(\rho)=\exp\left({1\over m}H_p(\rho)\right).
\end{eqnarray}

We need the following entropy dissipation formulas. In the case where $M$ is compact, the result is due to Lott and Villani. See  \cite{LV, V1, V2}

\begin{theorem} \label{AAAA} 
Let $(M, g)$ be a complete Riemannian manifold and $V\in C^2(M)$. 
Suppose that $Ric(L)=Ric+\nabla^2V$ is uniformly bounded on $M$. Let $d\mu=e^{-V}dv$, and
$$U (\rho)=\int_M e(\rho)d\mu,\ \ \rho\in P_2(M, \mu),$$ 
where $e\in C^2(\mathbb{R}^+, \mathbb{R}^+)$. Set\footnote{In Lott-Villani \cite{LV} and Villani \cite{V1, V2}, the function $p_1(r)$ was denoted by $p(r)$. Here we use the notation $p_1(r)$ instead of $p(r)$ to avoid the confusion of the exponent $p$ and the function $p(r)$.}  
\begin{eqnarray*}
p_1(r):=re'(r)-e(r), \ \ \ p_2(r):=rp_1'(r)-p_1(r). 
\end{eqnarray*}
Let $(\rho(t), \phi(t))$ be a smooth solution to the 
Benamou-Brenier geodesic equation on $P_2(M, e^{-V}dv)$ satisfying the following growth condition
\begin{eqnarray*}
 \int_M \left[|\nabla\phi|^2+|L\phi|^2+|e'(\rho)|^2+|\nabla e'(\rho)|^2+|L\phi|^2|p_1'(\rho)| \right]\rho \hspace{0.2mm} d\mu+\int_M {|\nabla e(\rho)^2\over \rho} \hspace{0.2mm} d\mu<\infty,
\end{eqnarray*}
and
\begin{eqnarray*}
 \int_M \left[|\nabla\phi|^2+|\nabla^2\phi|^2+|L\phi|^2+|\nabla L\phi|^2\right]p_1(\rho)\hspace{0.2mm} d\mu<\infty,
\end{eqnarray*} 
Then the following entropy dissipation formulas hold

\begin{eqnarray}
{d\over dt} U(\rho(t))
&=&\int_M  \nabla \phi \cdot  \nabla p_1(\rho)d\mu=-\int_M L\phi p_1(\rho) d\mu,\ \label{U2aa}
\end{eqnarray}
and
\begin{eqnarray}
{d^2\over dt^2} U(\rho(t))
=\int_M \Gamma_2(\phi, \phi) p_1(\rho)d\mu+\int_M (L \phi)^2 p_2(\rho)d\mu.  \label{U2a}
\end{eqnarray}where
\begin{eqnarray*}
\Gamma_2(\phi, \phi): ={1\over 2}L|\nabla \phi|^2- \nabla L \phi\cdot \nabla \phi =\|\nabla^2 \phi\|^2+Ric(L)(\nabla \phi, \nabla \phi).
\end{eqnarray*}
\end{theorem}
{\it Proof}. In the case where $M$ is compact, see  \cite{LV, V1, V2}. In the case where $(M, g)$ is a complete Riemannian manifold  and $(\rho, \phi)$ satisfies the required  growth condition at infinity, we can use the same argument as used for the proof of Theorem \ref{AAA} in \cite{LL24} to prove $(\ref{U2aa})$ 
and $(\ref{U2a})$. 

Let $\eta_k$ be an increasing sequence of functions in $C_0^\infty(M)$ such that $0\leq \eta_k\leq 1$, $\eta_k=1$ on $B(o, k)$ , $\eta_k=0$ on $M\setminus B(o, 2k)$, and $\|\nabla\eta_k\|\leq {1\over k}$. By standard argument and integration by parts, we have 
\begin{eqnarray*}
\partial_t \int_M e(\rho) \eta_k d\mu
&=&\int_M \partial_t e(\rho)\eta_k d\mu=\int_M e'(\rho)\partial_t\rho\eta_k d\mu\\
&=&\int_M \nabla_\mu^*(\rho\nabla\phi)e'(\rho)\eta_k d\mu=\int_M \langle \rho\nabla\phi, \nabla (e'(\rho)\eta_k)\rangle d\mu\\
&=&\int_M \rho \langle \nabla\phi, \nabla e'(\rho) \rangle \eta_k d\mu+\int_M \rho \langle\nabla\phi, e'(\rho)\nabla\eta_k\rangle d\mu\\
&:=&I_1(k)+I_2(k).
\end{eqnarray*}
Here
\begin{eqnarray*}
I_{1}(k)&=&\int_M \langle \nabla\phi, \nabla  e'(\rho) \rangle \eta_k \rho d\mu=
\int_M \langle \nabla\phi, \nabla  p_1(\rho) \rangle \eta_k  d\mu,\\
I_{2}(k)&=&\int_M\rho \langle\nabla\phi, e'(\rho)\nabla\eta_k\rangle d\mu.
\end{eqnarray*}

%
%
%

Under the assumption of theorem, we have $\int_M |\nabla\phi|^2\rho\hspace{0.2mm} d\mu<\infty$ and 
    $\int_M |\nabla e'(\rho)|^2 \rho\hspace{0.2mm} d\mu<\infty$, and hence $|\langle \nabla\phi, \nabla e'(\rho)\rangle| \in L^1(M, \rho\mu)$. 
By the Lebesgue dominated convergence theorem, we have
\begin{eqnarray}
I_1(k)\rightarrow \int_M \langle \nabla\phi, \nabla e'(\rho)\rangle \rho d\mu=\int_M \langle \nabla\phi, \nabla p_1(\rho)\rangle  d\mu. \label{en3}
\end{eqnarray}

Note that
\begin{eqnarray*}
I_{1}(k)&=&\int_M \nabla_\mu^*(\eta_k\rho\nabla\phi)e'(\rho)\hspace{0.2mm} d\mu\nonumber\\
&=&-\int_M \eta_k \rho L\phi  e'(\rho)\hspace{0.2mm} d\mu {\color{red}-} \int_M \langle \nabla \eta_k, \nabla\phi \rangle 
e'(\rho)\rho \hspace{0.2mm} d\mu {\color{red}-} \int_M  \eta_k\langle \nabla \rho, \nabla\phi \rangle e'(\rho)
\hspace{0.2mm}  d\mu.\nonumber\\
&=&-\int_M \eta_k \rho L\phi  e'(\rho)\hspace{0.2mm} d\mu {\color{red}-} \int_M \langle \nabla \eta_k, \nabla\phi \rangle 
e'(\rho)\rho \hspace{0.2mm} d\mu {\color{red}-} \int_M  \eta_k\langle \nabla e(\rho), \nabla\phi \rangle 
\hspace{0.2mm}  d\mu.
\end{eqnarray*}

Under the assumption of theorem, we have $\int_M |L\phi|^2\rho\hspace{0.2mm} d\mu<\infty$,  
$\int_M |e'(\phi)|^2\rho\hspace{0.2mm} d\mu<\infty$, and
 $\int_M {| \nabla e(\rho)|^2\over \rho}\hspace{0.2mm} d\mu<\infty$,
 $\int_M |\nabla\phi|^2\rho\hspace{0.2mm} d\mu<\infty$. Hence, the Cauchy-Schwarz inequality yields 
 $$\int_M |L\phi| |e'(\rho)|\rho\hspace{0.2mm} d\mu+\int_M |\nabla\phi| |e'(\rho)|\rho\hspace{0.2mm}
  d\mu+\int_M |\nabla e(\rho)||\nabla\phi| \hspace{0.2mm}  d\mu<\infty.$$

 Noting that $0\leq \eta_k\leq 1$, $\eta_k\rightarrow 1$ and 
$|\nabla\eta_k|\leq 1/k$, the Lebesgue dominated convergence theorem yields
\begin{eqnarray}
\lim\limits_{k\rightarrow \infty}I_1(k)
&=& -\int_M \rho L\phi e'(\rho)\hspace{0.2mm} d\mu {\color{blue}-} \int_M \langle \nabla e(\rho), \nabla\phi \rangle \hspace{0.2mm} d\mu\nonumber\\
&=&-\int_M \rho L\phi e'(\rho)\hspace{0.2mm} d\mu {\color{blue}+} \int_M  e(\rho) L\phi\hspace{0.2mm} d\mu \nonumber\\
&=&-\int_M  L\phi p_1(\rho)\hspace{0.2mm} d\mu. \label{en3b}
\end{eqnarray}
On the other hand, under the assumption of theorem, we have $\int_M |\nabla\phi|^2\rho\hspace{0.2mm} d\mu<\infty$ and $\int_M |e'(\rho)|^2\rho\hspace{0.2mm} d\mu<\infty$, which yields  $\int_M |e'(\rho)| |\nabla\phi| \rho\hspace{0.2mm} d\mu<\infty$.  By the Lebesgue dominated convergence theorem, as $|\nabla\eta_k|\leq 1/k$,  we have 
\begin{eqnarray}
I_2(k)=\int_M \rho\langle \nabla\phi, e'(\rho)\nabla \eta_k \rangle d\mu\rightarrow 0. \label{en4}
\end{eqnarray}
Combining $(\ref{en3})$,  $(\ref{en3b})$ with $(\ref{en4})$, we have

\begin{eqnarray*}
\partial_t U(\rho)&=&\lim\limits_{k\rightarrow \infty}[I_1(k)+I_2(k)]\\
&=&\int_M \langle \nabla \phi, \nabla p_1(\rho)\hspace{0.2mm} d\mu\\
&=&-\int_M L\phi p_1(\rho)\hspace{0.2mm} d\mu.
\end{eqnarray*}
This completes the proof of $(\ref{U2aa})$.  
%
%

Next, we prove $(\ref{U2a})$. By  standard argument, we have
\begin{eqnarray}
\partial_t \int_M L\phi p_1(\rho) \eta_k d\mu 
&=&\int_M \partial_t(L\phi p_1(\rho))\eta_kd\mu\nonumber\\
&=&\int_M L\partial_t \phi p_1(\rho) \eta_k d\mu+\int_M L\phi \partial_t p_1(\rho) \eta_kd\mu\nonumber\\
&=&-{1\over 2}\int_M L|\nabla\phi|^2p_1(\rho)\eta_k d\mu+\int_M L\phi p_1'(\rho) \nabla_\mu^*(\rho \nabla\phi)\eta_kd\mu\nonumber\\
&:=&I_3(k)+I_4(k). \label{en6}
\end{eqnarray}
Here
\begin{eqnarray*}
I_3(k)&=&-{1\over 2}\int_M L|\nabla\phi|^2 p_1(\rho) \eta_k d\mu,\\
I_4(k)&=&\int_M L\phi p_1'(\rho)\nabla_\mu^*(\rho \nabla\phi)\eta_kd\mu.
\end{eqnarray*}

By the weighted Bochner formula and $|Ric(L)|\leq C$, under the assumption $\int_M [|\nabla\phi|^2+|\nabla L\phi|^2+ |\nabla^2\phi |^2]p_1(\rho)\hspace{0.2mm} d\mu<\infty$, we have
\begin{eqnarray*}
\int_M |L|\nabla\phi|^2|p_1(\rho)\hspace{0.2mm} d\mu&=&2\int_M \left|\langle \nabla\phi, \nabla L\phi\rangle+|\nabla^2\phi|^2+Ric(L)(\nabla\phi, \nabla\phi)\right| p_1(\rho)\hspace{0.2mm} d\mu\\
&\leq&2 \int_M\left[|\nabla\phi| |\nabla L\phi|+|\nabla^2\phi|^2+C|\nabla\phi|^2\right| p_1(\rho)\hspace{0.2mm} d\mu<\infty.
\end{eqnarray*} 


Using  the fact $0\leq \eta_k\leq 1$ and $\eta_k\rightarrow 1$, the Lebesgue dominated convergence theorem yields
\begin{eqnarray}
I_3(k)\rightarrow -{1\over 2}\int_M L|\nabla\phi|^2p_1(\rho)\hspace{0.2mm} d\mu. \label{en7}
\end{eqnarray}

On the other hand, using the integration by parts formula, we have
\begin{eqnarray*}
I_4(k)
&=&\int_M L\phi p_1'(\rho) \nabla_\mu^*(\rho \nabla\phi)\eta_kd\mu\\
&=&\int_M L\phi p_1'(\rho) [-\rho L\phi-\langle\nabla \rho, \nabla \phi\rangle]\eta_k\hspace{0.2mm} d\mu\\
&=&-\int_M [|L\phi|^2 \rho p_1'(\rho)+L\phi \langle\nabla p_1(\rho), \nabla \phi\rangle]\eta_k\hspace{0.2mm} d\mu\\
&=&I_5(k)+I_6(k).
\end{eqnarray*}
Note that 
\begin{eqnarray*}
I_6(k)&=&-\int_M L\phi \langle\nabla p_1(\rho), \nabla \phi\rangle\eta_k\hspace{0.2mm} d\mu\\
&=&-\int_M \langle \eta_k L\phi \nabla \phi, \nabla p_1(\rho)\rangle\hspace{0.2mm} d\mu\\
&=&-\int_M  \nabla_\mu^*(\eta_k L\phi \nabla \phi) p_1(\rho)\hspace{0.2mm} d\mu\\
&=&\int_M [\eta_k |L\phi|^2+L\phi \nabla e_k \cdot\nabla \phi+\eta_k \nabla L\phi\cdot\nabla \phi]p_1(\rho)\hspace{0.2mm} d\mu.
\end{eqnarray*}

Under the assumption of theorem,  we have 
$$\int_M |L\phi|^2|p_1'(\rho)|\rho \hspace{0.2mm} d\mu+\int_M [ |L\phi|^2+|\nabla \phi|^2+ |\nabla L\phi|^2]p_1(\rho)\hspace{0.2mm} d\mu<\infty.$$

Using again the fact $0\leq \eta_k\leq 1$, $\eta_k\rightarrow 1$ and 
$|\nabla\eta_k|\leq 1/k$,  the Lebesgue dominated convergence theorem yields
\begin{eqnarray}
\lim\limits_{k\rightarrow \infty}I_4(k)&=&\lim\limits_{k\rightarrow \infty}[I_5(k)+I_6(k)]\nonumber\\
&=&-\int_M |L\phi|^2[\rho \langle p_1'(\rho)-p_1(\rho)]\hspace{0.2mm} d\mu+\int_M
\langle \nabla L\phi, \nabla\phi\rangle p_1(\rho)]\hspace{0.2mm} d\mu\nonumber\\
&=&-\int_M |L\phi|^2p_2(\rho)\hspace{0.2mm} d\mu+\int_M
\langle \nabla L\phi, \nabla\phi\rangle p_1(\rho)]\hspace{0.2mm} d\mu.\label{en8}
\end{eqnarray}

Combining $(\ref{en6})$, $(\ref{en7})$ and $(\ref{en8})$, we have

\begin{eqnarray*}
\partial_t^2 U(\rho)&=&-\lim\limits_{k\rightarrow \infty}[I_3(k)+I_4(k)]\\
&=&{1\over 2}\int_M L|\nabla\phi|^2p_1(\rho)\hspace{0.2mm} d\mu+\int_M |L\phi|^2 p_2(\rho)\hspace{0.2mm} d\mu -\int_M\langle \nabla L\phi, \nabla\phi\rangle p_1(\rho)\hspace{0.2mm} d\mu\\
&=&\int_M |L\phi|^2 p_2(\rho)\hspace{0.2mm} d\mu +\int_M \Gamma_2(\phi, \phi)p_1(\rho)\hspace{0.2mm} 
d\mu.
\end{eqnarray*}
This completes the proof of $(\ref{U2a})$. \hfill $\square$

%
%

\begin{theorem}\label{U2} Under the same condition as in Theorem \ref{AAAA}, we have
\begin{eqnarray}
{d^2\over dt^2} U(\rho(t))
&\geq &\int_M \left[\left(p_2(\rho)+{1\over m}p_1(\rho)\right)(L\phi)^2+Ric_{m, n}(L)(\nabla \phi, \nabla \phi) p_1(\rho)\right]d\mu.\nonumber\\
& &\label{UU2b1}
\end{eqnarray}
Assume that $p_2(r)+{1\over m}p_1(r)\geq \sigma p_1(r)$ and $Ric_{m, n}(L)\geq K$, then 
\begin{eqnarray}
{d^2\over dt^2} U(\rho(t))
\geq \sigma \left[{d\over dt}U(\rho(t))\right]^2\left[\int_M p_1(\rho(t))d\mu\right]^{-1}+\int_M K |\nabla \phi(t)|^2p_1(\rho(t))d\mu.\label{U3bb}
\end{eqnarray}
In particular,  if $p_2(r)+{1\over m}p_1(r)\geq 0$ and $Ric_{m, n}(L)\geq K$, we have
\begin{eqnarray}
{d^2\over dt^2} U(\rho(t))
\geq \int_M K |\nabla \phi(t)|^2p_1(\rho(t))d\mu.\label{U3bb0}
\end{eqnarray}
\end{theorem}
{\it Proof}. 
By $(\ref{U2a})$ and using the Bochner inequality 
\begin{eqnarray}
\Gamma_2(\phi, \phi)\geq {|L\phi|^2\over m}+Ric_{m, n}(L)(\nabla\phi, \nabla\phi),  \label{U2c}
\end{eqnarray}
we can derive $(\ref{UU2b1})$. Under the assumption $p_2(r)+{1\over m}p_1(r)\geq \sigma p_1(r)$ and the curvature-dimension condition $Ric_{m, n}(L)\geq Kg$, we have
\begin{eqnarray*}
{d^2\over dt^2} U(\rho(t))
\geq \sigma \int_M |L \phi|^2 p_1(\rho) d\mu+\int_M K |\nabla \phi|^2 p_1(\rho)d\mu. \label{U3a}
\end{eqnarray*}
By the Cauchy-Schwarz inequality
\begin{eqnarray*}
\left(\int_M L \phi  p_1(\rho) d\mu\right)^2\leq \left(\int_M |L \phi|^2  p_1(\rho) d\mu\right)\left(\int_M p_1(\rho)d\mu\right).
\end{eqnarray*}
Since $\sigma\geq 0$, we have
\begin{eqnarray*}
{d^2\over dt^2} U(\rho(t))
\geq \sigma \left(\int_M L \phi p_1(\rho) d\mu\right)^2\left(\int_M p_1(\rho)d\mu\right)^{-1}+\int_M K |\nabla \phi|^2 p_1(\rho)d\mu. \label{U3b}
\end{eqnarray*}
Using the fact $U'(\rho)=-\int_M p_1(\rho)L\phi d\mu$, we then obtain $(\ref{U3bb})$ \hfill $\square$

\subsection{Differential inequality for R\'enyi entropy}

In particular, for $e(r)={r^p\over p-1}$, we have  $U(\rho)={1\over p-1}\int_M \rho^p d\mu$, and
\begin{eqnarray*}
p_1(r)=re'(r)-e(r)=r^p, \ \ \ p_2(r)=rp'_1(r)-p_1(r)=(p-1)r^p.
\end{eqnarray*} 
Moreover
\begin{eqnarray*}
\nabla p_1(\rho)&=&\nabla \rho e'(\rho)+\rho e''(\rho)\nabla\rho-e'(\rho)\nabla\rho\\
&=&\rho e''(\rho)\nabla \rho.
\end{eqnarray*}
In this case, we have
\begin{eqnarray*}
p_2(\rho)+{1\over m}p_1(\rho)= \sigma p_1(\rho)
\end{eqnarray*}
with 
$$\sigma:=p-1+{1\over m}.$$

\begin{theorem}\label{Convex-3}  Let $e(r)={r^p\over p-1}$ and $U(\rho)={1\over p-1}\int_M \rho^p d\mu$. Assume that $Ric_{m, n}(L)\geq Kg$, where $K\in C(M, \mathbb{R})$. Then
\begin{eqnarray}
(p-1)U(\rho)U''(\rho)
\geq \left(p-1+{1\over m}\right) \left(U'(\rho)\right)^2+(p-1)U(\rho)\int_M K |\nabla \phi|^2 p_1(\rho)d\mu.\label{UP}
\end{eqnarray}
When $p\geq 1-{1\over m}$, we have
\begin{eqnarray}
{d^2\over dt^2} N_{m, p}(\rho)\leq -{N_{m, p}(\rho) \over m} {\int_M K |\nabla \phi|^2 \rho^p d\mu\over \int_M \rho^p d\mu}. \label{KNpa1}
\end{eqnarray}
If $Ric_{m, n}(L)\geq K$ for $K\in \mathbb{R}$, we have\begin{eqnarray*}
{d^2\over dt^2}N_{m, p}(\rho)\leq -{KN_{m, p}(\rho) \over m}
{ \int_M  |\nabla \phi|^2 \rho^p d\mu\over \int_M \rho^p d\mu}.
\end{eqnarray*}
In particular,  if $Ric_{m, n}(L)\geq 0$, i.e., if the CD$(0, m)$-condition holds, we have
\begin{eqnarray*}
{d^2\over dt^2}N_{m, p}(\rho)\leq 0.
\end{eqnarray*}
i.e., if $Ric_{m, n}(L)\geq 0$,  and $p\geq 1-{1\over m}$,  then $N_{m, p}(\rho(t))$ is concave in $t$ on $[0, \infty)$.  
\end{theorem}
{\it Proof}.  The inequality $(\ref{UP})$ follows from $(\ref{U3bb})$. Note that
\begin{eqnarray*}
H'_p(\rho)={1\over 1-p}{\partial_t \int_M \rho^pd\mu \over \int_M \rho^pd\mu}= {1\over 1-p}{U'\over U},
\end{eqnarray*}
and
\begin{eqnarray*}
H''_p(\rho)={1\over 1-p}{U''U-U'^2\over U^2}.
\end{eqnarray*}
Thus
\begin{eqnarray*}
H_p''+{1\over m}H_p'^2&=&{1\over p-1}{U'^2-U''U\over U^2}+{1\over m(p-1)^2}{U'^2\over U^2}\\
&=&{1\over (p-1)^2 U^2}\left[\left(p-1+{1\over m}\right)  U'^2- (p-1) U''U\right]\\
&\leq& -{ \int_M K |\nabla \phi|^2 p_1(\rho)d\mu\over \int_M p_1(\rho)d\mu}.
\end{eqnarray*}
Note that
\begin{eqnarray*}
{d\over dt}N_{m, p}(\rho)&=&{1\over m} H'_pN_{m, p}(\rho), \nonumber\label{NHp1}\\
{d^2\over dt^2}N_{m, p}(\rho)&=&{1\over m} \left(H_p''+{1\over m}H_p'^2\right)N_{m, p}(\rho). \label{NHp2}
\end{eqnarray*}
Combining the above calculations, we derive that
\begin{eqnarray*}
{d^2\over dt^2}N_{m, p}(\rho)\leq -{N_{m, p}(\rho) \over m}
{ \int_M K |\nabla \phi|^2 p_1(\rho)d\mu\over \int_M p_1(\rho)d\mu}.
\end{eqnarray*}
This proves $(\ref{KNpa1})$ since $p_1(\rho)=\rho^p$.  When $Ric_{m, n}(L)\geq K$ for $K\in \mathbb{R}$, we have\begin{eqnarray*}
{d^2\over dt^2}N_{m, p}(\rho)\leq -{KN_{m, p}(\rho) \over m}
{ \int_M  |\nabla \phi|^2 \rho^p d\mu\over \int_M \rho^p d\mu}.
\end{eqnarray*}
In particular,  if $Ric_{m, n}(L)\geq 0$, i.e., if the CD$(0, m)$-condition holds, we have
\begin{eqnarray*}
{d^2\over dt^2}N_{m, p}(\rho)\leq 0.
\end{eqnarray*}
\hfill $\square$

%
%
%
%
%
%
%
%
%
%
%
%

\medskip
\subsection{The $W$-entropy and NIW formula for R\'enyi entropy}

In this subsection, 
we introduce the $W$-entropy and prove the $W$-entropy formula for the R\'enyi entropy  along the geodesic flow on the Wasserstein space over complete 
Riemannian manifolds.
Then we prove the NIW formula which indicates the relationship between the entropy power $N$, the Fisher information $I=H'$ and the $W$-entropy related to the R\'enyi entropy for the geodesic
 flow on the Wasserstein space over 
compact Riemannian manifolds.  

When $m\in \mathbb{N}$, let $\rho_m(x, t)={e^{-{|x|^2\over 4t^2}}\over (2\pi t^2)^{m/2}}$, and $\phi_m(x, t)={|x|^2\over 2t}$. 
Then $(\rho_m, \phi_m)$ is a special solution to the geodesic equation on $P_2(\mathbb{R}^m, dx)$. See \cite{LL24}.
The R\'enyi entropyof $\rho_m(t)$ is given by 
\begin{eqnarray*}
H_p(\rho_m(t))&=&{1\over 1-p}\log \int_{\mathbb{R}^m} \rho_m(x, t)^p dx\\
&=&{1\over 1-p}\log \int_{\mathbb{R}^m} {e^{-{p|x|^2\over 4t^2}}\over (4\pi t^2)^{mp/2}}dx
\\
&=&{1\over 1-p} \left(\log \int_{\mathbb{R}^m} e^{-{p|x|^2\over 4t^2} }dx-
\log  (4\pi t^2)^{mp/2}\right)\\
\end{eqnarray*}
Using the fact 
 \begin{eqnarray*}
 \int_{\mathbb{R}^m} e^{-{p|x|^2\over 4t^2}}dx=(4\pi t^2p^{-1})^{m/2},
  \end{eqnarray*}
we have
  \begin{eqnarray*}
H_p(\rho_m(t))
={n\over 2}\log(4\pi et^2)+{m\over 2}{\log p\over p-1}-{m\over 2}.
\end{eqnarray*}
Note that
\begin{eqnarray*}
\lim\limits_{p\rightarrow 1}H_p(\rho_m(t))=H_1(\rho_m(t))={m\over 2}\log(4\pi et^2).
\end{eqnarray*}
%
%

\noindent{\bf Proof of Theorem \ref{W-entropy}}. To simplify the notation, we write $H_p'=H_p'(\rho(t))={d\over dt}H_p(\rho(t))$, etc. 
 Note that 
 \begin{eqnarray*}
 H_p'={1\over 1-p}{U'\over U}=\int_M L \phi d\gamma=I_p,
 \end{eqnarray*}  and
 \begin{eqnarray*}
 H_p''={1\over 1-p}  {U''U  -U'^2\over U^2}=-{U''\over \int_M p_1(\rho)d\mu}+ (p-1)\left (\int_M L \phi d\gamma\right)^2.
 \end{eqnarray*}
 Now
 \begin{eqnarray*}
U''&=& \int_M \left[\left( p-1+{1\over m }\right) |L \phi|^2+ Ric_{m, n}(L)(\nabla \phi, \nabla \phi)\right] p_1(\rho)d\mu\\
& & +\int_M\left( {m-n\over mn}\left(L \phi+{m\over m-n}\nabla V\cdot\nabla \phi\right)^2+\left\|\nabla^2 \phi-{\Delta \phi\over n} g\right\|^2\right)p_1(\rho)d\mu.
\end{eqnarray*}
 Therefore
 \begin{eqnarray*}
H''_p+{2\over t} H'_p-{m\over t^2} 
&=&- \int_M \left[\left( p-1+{1\over m }\right) |L \phi|^2+ Ric_{m, n}(L)(\nabla \phi, \nabla \phi)\right] d\gamma\\
& & -\int_M\left( {m-n\over mn}\left(L \phi+{m\over m-n}\nabla V\cdot\nabla \phi\right)^2+\left\|\nabla^2 \phi-{\Delta \phi\over n} g\right\|^2\right)d\gamma\\
& & + (p-1)\left (\int_M L \phi d\gamma\right)^2+{2\over t}\int_M L \phi d\gamma-{m\over t^2}.
\end{eqnarray*}
By 
\begin{eqnarray}
\left\|\nabla^2 \phi-{g\over t}\right\|_{\rm HS}^2
&=&{1\over n} \left|\Delta \phi-{n\over t}\right|^2+\left\|\nabla^2 \phi-{\Delta \phi\over n}g\right\|_{\rm HS}^2 \nonumber\\
&=&{1\over n} \left|L \phi-{m\over t}+\left(\nabla V\cdot\nabla \phi+{m-n\over t}\right)   \right|^2+\left\|\nabla^2 \phi-{\Delta \phi\over n}g\right\|_{\rm HS}^2  \nonumber\\
&=& {1\over m} \left|L \phi-{m\over t}  \right|^2-{1\over m-n}\left(\nabla V\cdot\nabla \phi+{m-n\over t}\right)^2 \nonumber \\
& &+{m-n\over mn}\left[L \phi+{m\over m-n}\nabla V\cdot\nabla \phi\right]^2+\left\|\nabla^2 \phi-{\Delta \phi\over n}g\right\|_{\rm HS}^2, \label{Hessv}
\end{eqnarray}
we have
 \begin{eqnarray*}
H''_p+{2\over t} H'_p-{m\over t^2}
&=&- \int_M \left[\left( p-1+{1\over m }\right) |L \phi|^2+ Ric_{m, n}(L)(\nabla \phi, \nabla \phi)\right] d\gamma\\
& & -\int_M \left\|\nabla^2 \phi-{g\over t}\right\|^2 d\gamma -{1\over m-n}\int_M \left(\nabla V\cdot\nabla \phi+{m-n\over t}\right)^2d\gamma\\
& & +{1\over m}\int_M \left|L \phi-{m\over t}\right|^2d\gamma+ (p-1)\left (\int_M L \phi d\gamma\right)^2+{2\over t}\int_M L \phi d\gamma-{m\over t^2}\\
&=&-(p-1) \left[\int_M  |L \phi|^2d\gamma-\left(\int_M L \phi d\gamma\right)^2\right]- \int_M  Ric_{m, n}(L)(\nabla \phi, \nabla \phi) d\gamma\\
& & -\int_M \left(\left\|\nabla^2 \phi-{g\over t}\right\|^2 +{1\over m-n} \left(\nabla V\cdot\nabla \phi+{m-n\over t}\right)^2\right)d\gamma.
\end{eqnarray*}

By the definition formula of $W_m$, we have 
\begin{eqnarray*}
{1\over t}{d\over dt}W_{m, p}(\rho, t)=H''_p+{2\over t} H'_p-{m\over t^2} .
\end{eqnarray*}
This completes the proof of $(\ref{derivativeofW})$. \hfill $\square$
%

\medskip

\begin{theorem} (NIW formula) Under the same notations as in Theorem \ref{Convex-3}, we have
\begin{eqnarray}
{m\over N_{m, p}(t)} {d^2\over dt^2} N_{m, p}(\rho(t))={1\over m}\left|I_p-{m\over t}\right|^2+{1\over t} {d\over dt}W_{m, p}(\rho, t).  \label{NIW-1}
\end{eqnarray}
In particular, when $Ric_{m, n}(L)\geq Kg$, we have
\begin{eqnarray}
 {d \over dt}W_{m, p}(\rho, t)\leq -t\int_M K|\nabla\phi|^2d\gamma-{t\over m} \left|I_p-{m\over t}\right|^2.  \label{NIW-2}
\end{eqnarray}
\end{theorem}
{\it Proof}. Note that
\begin{eqnarray*}
N_{m, p}''={1\over m} N_{m, p}\left(H_p''+{1\over m} H_p'^2\right),
\end{eqnarray*}
and 
\begin{eqnarray*}
H_p''+{1\over m}H_p'^2
&=&H_p''+{2\over t}H_p'-{m\over t^2}+{1\over m}H_p'^{2}-{2\over t}H_p'+{m\over t^2}\\
&=&{1\over t}{d\over dt}W_{m, p}(t)+{1\over m}\left|H_p'-{m\over t}\right|^2.
\end{eqnarray*}
This completes the proof of  $(\ref{NIW-1})$. By Theorem \ref{Convex-3} and $H_p'=I_p$, we obtain $(\ref{NIW-1})$.  \hfill $\square$

\subsection{Rigidity theorem of $W$-entropy for R\'enyi entropy}

Note that, when $m\in \mathbb{N}$ and $M=\mathbb{R}^m$, we have
\begin{eqnarray}
 {d^2\over dt^2}H_p(\rho_m(t))+{2\over t} {d\over dt}H_p(\rho_m(t))-{m\over t^2}=0.
\end{eqnarray}
and
\begin{eqnarray}
W_{m, p}(\rho, t)={d\over dt}(tH_{m, p}(\rho, t))=0.\label{Wmp0}
\end{eqnarray}

Due to this observation and inspired by Theorem \ref{W-Rigidity Shannon}, we have the following rigidity theorem of the $W$-entropy for the R\'enyi entropy along the Benamou-Brenier geodesic on the Wasserstein space over Riemannian manifold with CD$(0, m)$-condition. 

\begin{theorem} \label{W-Rigidity}  Let $(M, g, V)$ be a complete Riemannian manifold with bounded geometry condition as in Theorem \ref{AAAA}, and let $(\rho(t), \phi(t))$ be a smooth Benamou-Brenier geodesic on the Wasserstein space over Riemannian manifold satisfying  growth condition as required in Theorem \ref{AAAA}. Suppose that 
$Ric_{m, n}(L)\geq 0$, $p\geq 1$, and ${d\over dt}W_{m, p}(\rho, t)=0$ at some $t=\tau>0$. Then $M$ is isometric to $\mathbb{R}^n$, $m=n$, $V$ is a constant and  $(\rho, \phi)=(\rho_n, \phi_n)$ as given by \eqref{rigiditygeodesic}.
\end{theorem}
{\it Proof}. The proof is similar to the proof of Theorem \ref{W-Rigidity Shannon} given in \cite{LL24}. For the completeness of the paper, we reproduce it as follows. 

Under the condition $Ric_{m, n}(L)\geq 0$, and using the $W$-entropy formula $(\ref{derivativeofW})$, we see that 
${d\over dt}W_{m, p}(\rho, t)=0$ at some $t=\tau>0$ if and only if

 \begin{eqnarray*}
 {\rm Var}_\mu(L\phi)=0,\ \ Ric_{m, n}(L)(\nabla \phi, \nabla \phi)=0, 
 \end{eqnarray*}
 and
 \begin{eqnarray*}
 \left\|\nabla^2 \phi-{g\over \tau}\right\|^2 +{1\over m-n} \left(\nabla V\cdot\nabla \phi+{m-n\over \tau}\right)^2=0.
 \end{eqnarray*}
 Equivalently,  $\nabla^2\phi={g\over \tau}$, $Ric_{m, n}(L)(\nabla\phi, \nabla \phi)=0$ and $\nabla\phi\cdot\nabla V+{m-n\over \tau}=0$. 
By the Taylor expansion along geodesics,  $\phi (\cdot, \tau)$ has a unique minimum point $x_0\in M$, and every 
level set of $\phi$ is compact subset of $M$. Thus, $\phi (\cdot, \tau)$ is a strongly convex exhaustion function on 
 complete non-compact Riemannian manifold $M$. By Theorem 3 in \cite{GW}, $M$ is diffeomorphic to $\mathbb{R}^n$. Integrating $\nabla^2 \phi(\cdot, \tau)={g\over \tau}$ along the shortest geodesics between $x_0$ and 
$x$  on $M$ shows that
\begin{eqnarray*}
2\tau (\phi (x, \tau)-\phi (x_0, \tau))=d^2(x_0,x),\ \ \ \ \forall x\in M.
\end{eqnarray*} 
This yields 
\begin{eqnarray*}
\Delta d^2(x_0, x)=2n, \ \ \ \forall \ x\in M,
\end{eqnarray*}
which implies that $(M, g)$ is isometric to the Euclidean space $(\mathbb{R}^n, (\delta_{ij}))$. For this, see  \cite{Li12} and references therein. 
By the generalized Cheeger-Gromoll splitting theorem  (see Theorem 1.3, p. 565 in \cite{FLZ}), we can derive that $V$ must be a constant and $m=n$. 

Thus $\phi(\cdot, \tau)\in C^\infty(\mathbb{R}^n)$ satisfies $\nabla^2\phi(x, \tau)={\delta_{ij}\over \tau}$. This yields $\nabla\phi(x, \tau)={x\over \tau}$ under the assumption $\nabla\phi(0, \tau)=0$, and $\phi(x, \tau)={\|x\|^2\over 2\tau}$ up to an additional  constant.
By the Hopf-Lax formula for the solution to the Hamilton-Jacobi equation~\eqref{HJ} with $\phi(x, \tau)={\|x\|^2\over 2\tau}$, we have 
\begin{eqnarray*}
\phi(x, t)=\inf\limits_{y\in \mathbb{R}^n}\left\{{\|y\|^2\over 2\tau}+{\|x-y\|^2\over 2(t-\tau)}\right\}= {\|x\|^2\over 2t}, \ \ \ \ \ \forall \ t> \tau, \  x\in \mathbb{R}^n.
\end{eqnarray*}
By the uniqueness of the smooth solution  to the Hamilton-Jacobi equation~\eqref{HJ}, we see that $\phi(x, t)={\|x\|^2\over 2t}=\phi_n(x, t)$ for all $t>0$. 
Solving the transport equation~\eqref{TA} with the initial data $\lim\limits_{t\rightarrow 0}\rho(x, t)=\delta_0(x)$, we have 
$\rho(x, t) = {e^{-{\|x\|2\over 4t^2}}\over (4\pi t^2)^{n/2}}=\rho_n(x, t)$ for all $t>0$ and $x\in \mathbb{R}^n.$ \hfill $\square$

\subsection{Rigidity of R\'enyi entropy and entropy power}

In this subsection, using the NIW formula, we derive an explicit formula for the second derivative of the 
R\'enyi entropy and entropy power, and prove the rigidity theorems for the  R\'enyi entropy and entropy power inequality 
under the ${\rm CD}(K, m)$-condition. 
 
 \medskip

\noindent{\bf Proof of Theorem \ref{rigidity3}}. By $(\ref{U2a})$, we have
\begin{eqnarray*}
{d^2\over dt^2} U(\rho(t))=\int_M [|\nabla^2 \phi|^2+Ric(L)(\nabla \phi, \nabla \phi)] p_1(\rho)d\mu+\int_M (L \phi)^2p_2(\rho)d\mu.
\end{eqnarray*}
By the fact
\begin{eqnarray*}
\|\nabla^2 \phi\|^2&=&  {|L \phi|^2\over m}-{|\nabla V\cdot \nabla \phi|^2\over m-n}+\left\|\nabla^2 \phi-{\Delta \phi\over n} g\right\|^2\\
& &\hskip2cm +{m-n\over mn}\left(L \phi+{m\over m-n}\nabla V\cdot\nabla \phi\right)^2,
\end{eqnarray*}
we have

\begin{eqnarray}
{d^2\over dt^2} U(\rho(t))
&= &\int_M \left[{|L \phi|^2\over m}+Ric_{m, n}(L)(\nabla \phi, \nabla \phi)\right] p_1(\rho)d\mu+\int_M (L \phi)^2p_2(\rho)d\mu\nonumber\\
& & +\int_M\left( {m-n\over mn}\left(L \phi+{m\over m-n}\nabla V\cdot\nabla \phi\right)^2+\left\|\nabla^2 \phi-{\Delta \phi\over n} g\right\|^2\right)p_1(\rho)d\mu\nonumber\\
&=& \int_M \left( {1\over m }p_1(\rho)+p_2(\rho)\right) |L \phi|^2d\mu+\int_M Ric_{m, n}(L)(\nabla \phi, \nabla \phi) p_1(\rho)d\mu\nonumber\\
& & +\int_M\left( {m-n\over mn}\left(L \phi+{m\over m-n}\nabla V\cdot\nabla \phi\right)^2+\left\|\nabla^2 \phi-{\Delta \phi\over n} g\right\|^2\right)p_1(\rho)d\mu.\nonumber\\
& &\label{UU1}
\end{eqnarray}
Hence,  for $e(r)={r^p\over p-1}$, we have
\begin{eqnarray}
{d^2\over dt^2} U(\rho(t))
&=& \int_M \left[\left( p-1+{1\over m }\right) |L \phi|^2+ Ric_{m, n}(L)(\nabla \phi, \nabla \phi)\right] p_1(\rho)d\mu\nonumber\\
& & +\int_M\left( {m-n\over mn}\left(L \phi+{m\over m-n}\nabla V\cdot\nabla \phi\right)^2+\left\|\nabla^2 \phi-{\Delta \phi\over n} g\right\|^2\right)p_1(\rho)d\mu.\nonumber\\
& & \label{UU2}
\end{eqnarray}
This yields
\begin{eqnarray*}
H_p''+{1\over m} H_p'^2
&=&{1\over (p-1)^2 U^2}\left[\left(p-1+{1\over m} \right)  U'^2- (p-1) U''U\right]\\
&=&\left(p-1+{1\over m} \right)\left({U'\over (p-1)U}\right)^2-{U''\over (p-1)U}\\
&=& \left(  p-1+{1\over m} \right) \left[\left(  \int_M L \phi d\gamma \right)^2-   \int_M  |L \phi|^2 d\gamma \right]-\int_M Ric_{m, n}(L)(\nabla \phi, \nabla \phi) d\gamma\\
& &\hskip1cm -\int_M \left[ {m-n\over mn}\left(L \phi+{m\over m-n}\nabla V\cdot\nabla \phi\right)^2 + \left\|\nabla^2 \phi-{\Delta \phi\over n} g\right\|^2_{\rm HS} \right]d\gamma.
\end{eqnarray*}
By the fact 
$${d^2\over dt^2}N_{m, p}(\rho(t))={1\over m}\left(H_p''+{1\over m}H_p'^2\right)N_{m, p}(\rho(t)),$$
we  can derive $(\ref{Np2})$.  
{By the same argument as used in the proof of the rigidity part of Theorem \ref{thm2}, the rigidity part  of Theorem \ref{rigidity3}  can be easily derived from $(\ref{Np2})$. } \hfill $\square$

In particular, we have the following rigidity theorem.

\begin{theorem} \label{rigidity4} If $Ric_{m, n}(L)\geq 0$ and $p\geq 1-{1\over m}$, we have the enhanced entropy power differential inequality
\begin{eqnarray}
& &H_p''+{1\over m} H_p'^2+\left(  p-1+{1\over m} \right) {\rm Var}_\gamma(L \phi)\nonumber\\
& &\hskip0.5cm +\int_M \left[ {m-n\over mn}\left(L \phi+{m\over m-n}\nabla V\cdot\nabla \phi\right)^2 + \left\|\nabla^2 \phi-{\Delta \phi\over n} g\right\|^2_{\rm HS} \right]d\gamma\nonumber\\ 
& &\hskip3cm \leq 0, \label{Hpm2K=0}
\end{eqnarray}  and
\begin{eqnarray}
& &{d^2\over dt^2}N_{m, p}(\rho(t))+\left(  p-1+{1\over m} \right) { N_{m, p}(\rho(t))\over m}{\rm Var}_\gamma(L \phi)\nonumber\\
& &\hskip1cm +{N_{m, p}(\rho(t))\over m} \int_M \left[ {m-n\over mn}\left(L \phi+{m\over m-n}\nabla V\cdot\nabla \phi\right)^2 + \left\|\nabla^2 \phi-{\Delta \phi\over n} g\right\|^2_{\rm HS} \right]d\gamma\nonumber\\
& &\hskip4cm \leq 0,
 \label{Np2K=0}
\end{eqnarray}
and $N_{m, p}(\rho(t))$ is concave in $t$ on $[0, \infty)$, i.e., 
\begin{eqnarray}
{d^2\over dt^2} N_{m, p}(\rho(t))\leq 0, \label{KNp0}
\end{eqnarray}
{
Moreover, the equality in the enhanced entropy differential inequality $(\ref{Hpm2K=0})$ or 
the enhanced entropy power differential inequality $(\ref{Np2K=0})$ holds for all geodesics $(\rho, \phi)$ at some $t$ if and only if $(M, g, V)$ is $(0, m)$-Einstein, i.e., 
$$Ric_{m, n}(L)=0.$$ 
{Futhermore, under the condition $Ric_{m, n}(L)=0$ and $p\geq 1-{1\over m}$,  the equality  in
$(\ref{Hpm2K=0})$ or  $(\ref{KNp0})$ holds at some $t$ if and only if $\phi$ satisfies the Hessian soliton equation
\begin{eqnarray*}
L\phi=I_p=\int_M L\phi d\gamma, \ \ L \phi+{m\over m-n}\nabla V\cdot\nabla \phi=0, \ \ \nabla^2 \phi={I_p\over m}g.
\end{eqnarray*}
In the case $(M, g, V)$ is a complete Riemannian manifold with bounded geometry condition as in Theorem \ref{AAAA}, and $(\rho(t), \phi(t))$ is a smooth Benamou-Brenier geodesic on the Wasserstein space over Riemannian manifold satisfying  growth condition as required in Theorem \ref{AAAA}. Then, assumption $Ric_{m, n}(L)=0$ and $p\geq 1-{1\over m}$,  the equality  in $(\ref{Hpm2K=0})$ or $(\ref{KNp0})$ holds at some $t=\tau>0$ if and only if $M$ is isometric to $\mathbb{R}^n$, $m=n$, $V$ is a constant and  $(\rho, \phi)=(\rho_n, \phi_n)$ as given by \eqref{rigiditygeodesic}.}

}
\end{theorem}
{\it Proof}. 
{The last statement of rigidity theorem is indeed a reformulation of Theorem 
\ref{W-Rigidity}.} \hfill $\square$

\medskip
In particular, for $m=n$, $V=0$, $\mu=v$, $L=\Delta$ and $K=0$, we have  
\begin{theorem} Under the assumption $Ric\geq 0$ and $p\geq 1-{1\over n}$,  we have the enhanced entropy power differential inequality
we have
\begin{eqnarray}
& &{d^2\over dt^2}N_{n, p}(\rho(t))+ \left(  p-1+{1\over n} \right)  {N_{n, p}(\rho(t))\over n}{\rm Var}_\gamma(\Delta \phi)+ {N_{n, p}(\rho(t))\over n} \int_M \left\|\nabla^2 \phi-{\Delta \phi\over n} g\right\|^2_{\rm HS}d\gamma\nonumber\\
& &\hskip4cm \leq 0. \label{NIWK=0}
\end{eqnarray}
In particular, 
\begin{eqnarray}
{d^2\over dt^2} N_{n, p}(\rho(t))\leq 0,\label{NpK0}
\end{eqnarray}
and
\begin{eqnarray}
 {d \over dt}W_{n, p}(\rho, t)\leq -{t\over n} \left|I_p-{n\over t}\right|^2. \label{NIW-2n0}
\end{eqnarray}
As a consequence,  $N_{n, p}(\rho(t))$ is concave in $t$ on $[0, \infty)$, and $W_{n, p}(\rho, t)$ is nonincreasing in $t$ on $[0, \infty)$. 

\medskip
Moreover,  the equality in the enhanced entropy power differential inequality $(\ref{NIWK=0})$ holds for all geodesics $(\rho, \phi)$ at some $t$ if and only if $(M, g)$ is Ricci flat in the sense that
$$Ric=0.$$ 
Furthermore, under the condition $Ric=0$, \eqref{NpK0} or \eqref{NIW-2n0} become an equality if and only if $\phi$ satisfies the Hessian soliton equation
\begin{eqnarray*}
\Delta \phi=I_p:=\int_M \Delta \phi d\gamma, \ \ \nabla^2 \phi={I_p\over n}g.
\end{eqnarray*}
Suppose that $(M, g)$ is a complete Riemannian manifold with bounded geometry condition as in Theorem \ref{AAAA}, and $(\rho(t), \phi(t))$ is a smooth Benamou-Brenier geodesic on the Wasserstein space over Riemannian manifold satisfying  growth condition as required in Theorem \ref{AAAA}. Then, under the assumption $Ric=0$ and $p\geq 1-{1\over n}$,  the equality  in \eqref{NpK0} or  $(\ref{NpK0})$ holds at some $t=\tau>0$ if and only if $M$ is isometric to $\mathbb{R}^n$ and  $(\rho, \phi)=(\rho_n, \phi_n)$ as given by \eqref{rigiditygeodesic}.

\end{theorem}

\subsection{The special case $p=1-1/N$ as in Sturm \cite{Sturm06b}}

In the special case  $e(r)={r^p\over p-1}$ with $p=1-{1/N}$, we have 
$$\sigma:=p-1+{1\over m }={1\over m}-{1\over N}\geq 0$$ if $N\geq m$. 
Using the same notation as in Sturm \cite{Sturm06b}, let 
$$S_{N}=-\int_M \rho^{1-1/N}d\mu.$$ Then 
$$U(\rho)={1\over p-1}\int_M \rho^pd\mu=NS_{N}(\rho).$$
Theorem \ref{Convex-3} and $(\ref{UU2})$ imply the following

\begin{corollary} Under the above notation, for all $N\geq m$, we have
\begin{eqnarray*}
N{d^2\over dt^2}S_N(\rho)
&=& \int_M \left[\left({1\over m }-{1\over N}\right) |L \phi|^2+ Ric_{m, n}(L)(\nabla \phi, \nabla \phi)\right] 
\rho^{1-1/N}d\mu\\
& & +\int_M\left[ \left({1\over n}-{1\over m}\right)\left(L \phi+{m\over m-n}\nabla V\cdot\nabla \phi\right)^2+\left\|\nabla^2 \phi-{\Delta \phi\over n} g\right\|^2\right] \rho^{1-1/N}d\mu.
\end{eqnarray*}
Assume that $Ric_{m, n}(L)\geq K$ for some $K\in C(M, \mathbb{R})$. Then, for all $N\geq m$, we have
\begin{eqnarray*}
{d^2\over dt^2}S_N(\rho(t))+{N-m\over m}[S_N(\rho(t))]^{-1}\left({d\over dt}S_N(\rho(t))\right)^2
\geq {1\over N} \int_M K |\nabla \phi|^2 \rho^{1-1/N}d\mu.
\end{eqnarray*}
In particular, for all $N\geq m$, we have
\begin{eqnarray}
{d^2\over dt^2}S_N(\rho(t))\geq {1\over N}\int_M K |\nabla \phi|^2 \rho^{1-1/N}d\mu.\label{SNK1a}
\end{eqnarray}
\end{corollary}

We will see in Section 5 that $(\ref{SNK1a})$ is the differential form of  $(\ref{S(K, N)})$. Hence 
 $(vi)$ in Theorem \ref{Thm2} implies $(viii)$.

\subsection{Limiting case $N\rightarrow \infty$ and $p\rightarrow 1$}

Letting $N\rightarrow \infty$ in $(\ref{SNK1a})$, and using
 \begin{eqnarray*}
\lim\limits_{N\rightarrow \infty} N [1+S_{N}(\rho)]=\int_M \rho \log\rho d\mu,
\end{eqnarray*}
we can derive 

 \begin{eqnarray*}
{d^2\over dt^2}\int_M \rho d\rho d\mu\geq K\int_M |\nabla \phi|^2 d\mu.\label{SN=1}
\end{eqnarray*}
i.e., 
 \begin{eqnarray*}
{d^2\over dt^2}{\rm Ent}(\rho(t))\geq KW_2^2(\rho(0), \rho(1)),
\end{eqnarray*}
That is to say, under CD$(K, \infty)$-condition, we have
 \begin{eqnarray*}
{\rm Hess}~{\rm Ent}(\rho)\geq K.
\end{eqnarray*}

\section{Comparison with Lott-Villani-Sturm's synthetic geometry}

%
%
The aim of this section is to compare the main results of this paper with the works of 
Lott-Villani \cite{LV}, Sturm \cite{Sturm06a, Sturm06b} and Villani \cite{V2} in the setting of Riemannian manifolds with $N=m$.

%
\subsection{Lott-Villani and Sturm}

%
%
For  $N<\infty$,  Sturm \cite{Sturm06b} proved that in the Riemannian case the curvature-dimension CD$(0, N)$ for 
$N\geq n={\rm dim} M$ holds if and only if for all $N'\geq N$ the quantity $S_{N'}(\cdot|m)$ is convex 
on the Wasserstein space $P_2(M, d)$.  This is basically due to the fact that the Jacobian determinant $J_t={\rm det}dF_t$ of any optimal transport map 
$F_t=\exp(-t\nabla\phi): M\rightarrow M$ satisfies 
\begin{eqnarray}
{\partial^2\over \partial t^2}J_t^{1/N}\leq 0\label{JN0}
\end{eqnarray}
if and only if $N\geq n={\rm dim} M$ and $Ric\geq 0$. As pointed out in \cite{Sturm06b}, $(\ref{JN0})$ is essentially equivalent to the Brunn-Minkowski inequality 
\begin{eqnarray}
m(A_t)^{1/N'}\geq (1-t)m(A_0)^{1/N'}+tm(A_1)^{1/N'}.  \label{BM1}
\end{eqnarray}
for any $N'\geq N$, any $t\in [0, 1]$ and any pair of Borel sets $A_0, A_1\subset M$, where $A_t$ denotes the set of 
points $\gamma_t$ on the geodesics with endpoints $\gamma_0\in A_0$ and $\gamma_1\in A_1$.
%
%

The curvature-dimension condition CD$(K, N)$  for general $K\in \mathbb{R}$ and $N<\infty$ is more involved. In this case, the inequality $(\ref{JN0})$ is replaced by 
\begin{eqnarray}
{\partial^2\over \partial t^2}J_t^{1/N}(x)\leq -{K\over N}J_t^{1/N}(x)
d^2(x, F_1(x)).
\label{JNK}
\end{eqnarray}
When $\mu=v$ is the standard volume measure, Sturm \cite{Sturm06b} proved that $(\ref{JNK})$ is equivalent to $(\ref{S(K, N)})$ if and only if $N\geq n={\rm dim}$ and $Ric\geq K$. In \cite{Sturm06b}, Sturm used $(\ref{S(K, N)})$  to 
 introduce the definition of the CD$(K, N)$-condition on metric measure spaces. 
Independent of Sturm \cite{Sturm06b}, Lott and Villani \cite{LV}  introduced the curvature-dimension condition CD$(0,N)$ on metric measure spaces in the same form as  $(\ref{S(0, N)})$.

By Sturm \cite{Sturm06b}, assume for simplicity that for 
$m\otimes m$ -a.e. $(x, y)\in M^2$ there exists a unique 
geodesic$t\rightarrow \gamma_t(x, y)$ depending in a 
measurable way on the endpoints $x$ and $y$. Then CD$(K,N)$
states that that for any pair of absolutely continuous probability 
measures  $\rho_0 m$ and  $\rho_1 m$
on $M$ there exists an optimal coupling $q$ such that

\begin{eqnarray}
\rho_t^{-1/N}(\gamma_t(x, y))\leq \tau_{K, N}^{(1-t)}(d(x, y))
\rho_0^{-1/N}(x)+\tau_{K, N}^{(t)}(d(x, y))
\rho_1^{-1/N}(y) \label{rho1}
\end{eqnarray}
for all $t\in [0, 1]$ and $q$-a.e. $(x, y)\in M^2$, where $\rho_t$ 
is the density of the push-forward of $q$ under
the map $(x, y)\rightarrow\gamma_t(x, y)$. 
%
Indeed, $\rho_t m=(F_t)_{*}(\rho_0 m)$. More precisely, we have

\begin{eqnarray*}
\rho_0(x)=\rho_t(F_t(x)){\rm det}(dF_t(x))=\rho_t(F_t(x))J_t. \label{rho2}
\end{eqnarray*}
Equivalently

\begin{eqnarray*}
\rho_t(F_t(x)) =\rho_0(x)J_t^{-1}. \label{rho3}
\end{eqnarray*}
Hence
\begin{eqnarray*}
{d^2\over dt^2}\rho_t^{-1/N}(F_t(x)) &=&\rho_0^{-1/N}(x){d^2\over dt^2}J_t^{1/N}\\
&\leq& {Kd^2(x, F_1(x))\over N}\rho_0^{-1/N}(x)J_t^{1/N}.
\end{eqnarray*}
That is to say
\begin{eqnarray}
{d^2\over dt^2}\rho_t^{-1/N}(F_t(x))\leq {Kd^2(x, F_1(x))\over N}\rho_t^{-1/N}(F_t(x)). \label{rho4}
\end{eqnarray}
Integrating in $t$, $(\ref{rho4})$ implies $(\ref{rho1})$.

%
%

Moreover, one can show that $(\ref{rho4})$ is equivalent to $(\ref{S(K, N)})$. Indeed, if  $(\ref{rho4})$ hold, then

\begin{eqnarray*}
{d^2\over dt^2}S_N(\rho(t))
&=&-{d^2\over dt^2}\int_M \rho_t^{-1/N}(F_t(x))  \rho_0(x)d\mu(x)\\
&=&-\int_M {d^2\over dt^2}\rho_t^{-1/N}(F_t(x))  \rho_0(x)d\mu(x)\\
&\leq&-{1\over N}\int_M K(x)d^2(x, F_1(x)) \rho_t^{-1/N}(F_t(x))  \rho_0(x)d\mu(x).
\end{eqnarray*}
Hence

\begin{eqnarray*}
{d^2\over dt^2}S_N(\rho(t))\leq-{1\over N}\int_M K(x)d^2(x, F_1(x))\rho_t^{1-1/N}(F_t(x)) J_t(x)d\mu(x).
\end{eqnarray*}
When $K$ is a  constant, we have 

\begin{eqnarray}
{d^2\over dt^2}S_N(\rho(t))\leq-{K\over N}\int_M d^2(x, F_1(x))\rho_t^{1-1/N}(F_t(x))J_t(x) d\mu(x).\label{SNK1}
\end{eqnarray}
On the other hand, if $(\ref{SNK1})$ holds for all $\rho_0(x)$, we can derive that $(\ref{rho4})$ also holds. 

Now  we prove that 

\begin{eqnarray}
\int_M |\nabla \phi_t(x)|^2\rho_t^{1-1/N}(x) d\mu(x)
=\int_M d^2(x, F_1(x))\rho_t^{1-1/N}(F_t(x))J_t(x) d\mu(x).\label{IJ}
\end{eqnarray}
Indeed

\begin{eqnarray*}
\int_M |\nabla \phi_t(x)|^2\rho_t^{1-1/N}(x) d\mu(x)&=&\int_M |\nabla \phi_t(y)|^2\rho_t^{1-1/N}(y) d\mu(y)\\
&=&\int_M |\nabla \phi_t(F_t(x))|^2\rho_t^{1-1/N}(F_t(x)){\rm det}(dF_t(x))d\mu(x)\\
&=&\int_M |\nabla \phi_t(F_t(x))|^2\rho_t^{-1/N}(F_t(x))\rho_0(x)d\mu(x).
\end{eqnarray*}
Note that $\nabla \phi_t(F_t(x))={d\over dt}F_t(x)=\dot F_t(x)$, and $F_t(x)=\exp_x\left(-t\nabla \phi(x)\right)$ is a geodesic linking $F_0(x)=x$ and $F_1(x)$ . Hence 
\begin{eqnarray*}
|\nabla\phi_t(F_t)|^2=|\dot F_t(x)|^2=|\dot F_1(x)|^2=d^2(x, F_1(x)).
\end{eqnarray*}
This proves $(\ref{IJ})$. As a consequence, $(\ref{SNm4b})$) for $S_N(\rho(t))$ is the same as  $(\ref{SNK1})$, which is equivalent to Sturm's definition inequality $(\ref{S(K, N)})$. This proves $(vii)$ in Theorem \ref{Thm2} implies $(viii)$ there. 

Next we briefly describe the works of  Lott-Villani \cite{LV} and Villani \cite{V2} on the characterization of CD$(K, N)$-condition by $\mathcal{DC}_N$ class of functional.  In \cite{LV, V2}, Lott and Villani introduced the class $\mathcal{DC}_N$ of general functional $U(\rho)=\int_M e(\rho)d\mu$ on the Wasserstein space $P_2(M)$. By Definition 17.1 in \cite{V2}, the class $\mathcal{DC}_N$ is defined as 
the set of continuous convex function  $e:\mathbb{R}^+\rightarrow 
\mathbb{R}$ which is $C^2$-differentiable on $(0, \infty)$, $e(0)=0$ and satisfies one of the following equivalent conditions holds:
\begin{itemize}

\item The pressure function $p_1(r)=re'(r)-e(r) $ satisfies
 \begin{eqnarray*}
r p_1'(r)\geq \left(1-{1\over N}\right)p_1(r). \label{p1}
\end{eqnarray*}
Equivalently
\begin{eqnarray*}
p_2(r)+{1\over N}p_1(r)\geq 0. \label{p2}
\end{eqnarray*}

\item The function $r\mapsto {p_1(r)\over r^{1-{1\over N}}}$ is convex. 

\item The function $\delta\rightarrow u(\delta)=\delta^{N}e(\delta^{-N})$ is convex. 

\end{itemize}

By Theorem 17.15  in \cite{V2}, $M$ satisfies the curvature-dimension condition CD$(K, N)$ if and only if for each 
$e\in \mathcal{DC}_N$, the functional $U(\rho)=\int_M e(\rho)d\mu$ satisfies
\begin{eqnarray}
U(\mu_t)+K_{N, U}\int_0^1 \int_M \rho_s(x)^{1-{1\over N}}|\nabla \phi_s(x)|^2d\mu(x)G(s, t)ds \leq (1-t)U(\mu_0)+tU(\mu_1). \label{Ut1}
\end{eqnarray} 
Here
$$K_{N, U}=\inf\limits_{r>0} {Kp_1(r)\over r^{1-1/N}}.$$
Indeed, it was proved in \cite{V2} that for any $e\in \mathcal{DC}_N$, it holds
\begin{eqnarray}
{d^2\over dt^2}U(\mu_t)\geq K_{N, U}\int_M |\nabla \phi_t(x)|^2\rho_t(x)^{1-1/N}d\mu(x).\label{Ut2}
\end{eqnarray} 
which is equivalent to $(\ref{Ut1})$ by Proposition 16.2 in Villani \cite{V2}.

Now we compare Theorem \ref{U2} and Lott-Villani's inequalities $(\ref{Ut1})$ and $(\ref{Ut2})$. 
Note that, by $(\ref{U3bb})$ in Theorem \ref{U2}, if $p_2+{1\over N}p_1\geq \sigma p_1$ for some constant 
$\sigma\geq 0$, it holds
   \begin{eqnarray*}
 U''(\rho_t)\geq \sigma [U'(\rho_t)]^2 \left(\int_M p_1(\rho_t(x))d\mu(x)\right)^{-1}+ \int_M K(x)|\nabla \phi_t(x)|^2p_1( \rho_t(x))d\mu(x). \label{Uttt}
  \end{eqnarray*}
  In particular, for $e\in \mathcal{DC}_N$,  we have $\sigma=0$. Hence 
   \begin{eqnarray*}
 U''(\rho_t)\geq \int_M K(x)|\nabla \phi_t(x)|^2p_1( \rho_t(x))d\mu(x).\label{Uttt2}
  \end{eqnarray*}
 Under Lott-Villani's assumption $Kp_1(\rho)\geq K_{N, U}\rho^{1-1/N}$, 
$(\ref{Ut2})$ follows.

\begin{remark}
Similarly to Lott-Villani \cite{LV} and Villani \cite{V2}, we may introduce the class 
$\mathcal{DC}_{N, \sigma}$ of functionals $U(\rho)=\int_M e(\rho)d\mu$ such that 
$p_2(r)+{1\over N}p_1(r)\geq \sigma p_1(r)$ for a constant $\sigma\geq 0$. Note that $\mathcal{DC}_{N, \sigma_2}\subset 
\mathcal{DC}_{N, \sigma_1}$ if $\sigma_2\geq \sigma_1\geq 0$, and 
$e(r)={r^p\over p-1}\in \mathcal{DC}_{N, \sigma}$ with
$\sigma=p-1+{1\over N}$. 
\end{remark}
%
%
%

%
%
%

\subsection{Comparison}

Now we compare our entropy power differential inequalities  in Theorem \ref{Thm2} with the above mentioned results due to Sturm \cite{Sturm06b}, Lott-Villani \cite{LV} and Villani \cite{V2}. 

Note that our entropy power inequality  
 \begin{eqnarray}
 {d^2\over dt^2}N_{m}(\rho(t))\leq -{K\over m}N_m(\rho(t))W_2^2(\rho(0), \rho(1))  \label{Ncon2}
 \end{eqnarray}
is equivalent to the entropy differential inequality
 \begin{eqnarray}
 H''+{H'^2\over m}\leq K\int_M |\nabla\phi|^2\rho d\mu, \label{Hm0}
 \end{eqnarray} 
 where $(\rho(t), t\in [0, 1])$ is a geodesic linking $\rho(0)$ and $\rho(1)$ on the Wasserstein space $P_2(M, g, \mu)$, $H'={d\over dt}H(\rho(t))$ and $H''={d^2\over dt^2}H(\rho(t))$. 
 
From the proof of the NIW formula, we have
 
\begin{eqnarray}
H''+{H'^2\over m}
&=&-{1\over m} \int_M \left| L \phi- \int_M L \phi  \rho d\mu \right|^2 \rho d\mu-\int_M Ric_{m, n}(L)(\nabla \phi, \nabla \phi) \rho d\mu\nonumber\\
& &-{m-n\over mn} \int_M \left(L \phi+{m\over m-n}\nabla V\cdot\nabla \phi\right)^2\rho d\mu- \int_M \left\|\nabla^2 \phi-{\Delta \phi\over n} g\right\|^2_{\rm HS}\rho d\mu. \nonumber\\
& &\label{hhhh0}
\end{eqnarray}
Under the assumption $Ric_{m, n}(L)\geq K$, we have 
\begin{eqnarray}
& &H''+{H'^2\over m}+{1\over m} \int_M \left| L \phi- \int_M L \phi  \rho d\mu \right|^2 \rho d\mu\nonumber\\
& &\hskip1cm+{m-n\over mn} \int_M \left(L \phi+{m\over m-n}\nabla V\cdot\nabla \phi\right)^2\rho d\mu- \int_M \left\|\nabla^2 \phi-{\Delta \phi\over n} g\right\|^2_{\rm HS}\rho d\mu. \nonumber\\
& &\hskip2cm \leq -K\int_M |\nabla \phi|^2\rho d\mu=-KW_2^2(\rho_0, \rho_1),\label{hhhh0H}
\end{eqnarray}
which yields
\begin{eqnarray}
N_m''\leq -{KN_m\over m}W_2^2(\rho_0, \rho_1).
\label{hhhh0N}
\end{eqnarray}

In particular, when $m=n$, it holds

\begin{eqnarray}
H''+{H'^2\over n}
&=&-{1\over n} \int_M \left| \Delta \phi- \int_M \Delta \phi  \rho d\mu \right|^2 \rho d\nu-\int_M Ric(\nabla \phi, \nabla \phi) \rho d\nu \nonumber\\\
& &\hskip2cm - \int_M \left\|\nabla^2 \phi-{\Delta \phi\over n} g\right\|^2_{\rm HS}\rho d\nu\nonumber\\\
&=&-\int_M \left[Ric(\nabla \phi, \nabla \phi) +\left\|\nabla^2 \phi-{\int_M \Delta \phi \rho d\nu\over n} g\right\|^2_{\rm HS}\right]\rho d\nu.  \label{hhhh1}
\end{eqnarray}
Hence, under the condition $Ric\geq K$, we have 
\begin{eqnarray}
& &H''+{H'^2\over n}+{1\over n} \int_M \left| \Delta \phi- \int_M \Delta \phi  \rho d\mu \right|^2 \rho d\nu+\int_M \left\|\nabla^2 \phi-{\int_M \Delta \phi \rho d\nu\over n} g\right\|^2_{\rm HS}\rho d\nu \nonumber\\
& &\hskip2cm \leq -K \int_M |\nabla \phi|^2\rho d\nu=-KW_2^2(\rho(0), \rho(1)),\label{hhhh2}
\end{eqnarray}
which yields
\begin{eqnarray}
N_n''\leq -{KN_m\over n}W^2(\rho_0, \rho_1).
\label{hhhh2N}
\end{eqnarray}

%
%
%
Conversely, we now prove that,  if $(\ref{hhhh2})$ or $(\ref{hhhh2N})$ holds for any geodesic on the Wasserstein 
space $P_2(M, v)$ over $(M, g)$, then we have 
\begin{eqnarray*}
Ric\geq K.
\end{eqnarray*}
Indeed, $(\ref{hhhh1})$ holds at $t=0$ for the solution $(\rho(t), \phi(t))$ of the continuity equation and the Hamilton-Jacobi equation with any initial data 
$(\rho(0), \phi(0))$. Taking a sequence of $\rho_k(0)d\nu$ converges weakly to the Dirac mass at $x_0\in M$, we derive that at $t=0$

\begin{eqnarray*}
H''+{H'^2\over n}=-\left[Ric_{x_0}(v, v)+\left\|\nabla^2 \phi(x_0)- {\Delta \phi(x_0)\over n} g\right\|^2_{\rm HS}\right].
\end{eqnarray*}
For any fixed $v\in T_{x_0}M$, choose  $\phi\in C^2(M)$ such that 
\begin{eqnarray*}
\nabla \phi(x_0)&=&v,\\
\nabla^2\phi(x_0)&=&{\Delta \phi(x_0)\over n}I_n.
\end{eqnarray*}
This implies that

\begin{eqnarray*}
H''+{H'^2\over n}=-Ric_{x_0}(v, v).
\end{eqnarray*}
Therefore, if $(\ref{hhhh2})$ holds, we have  

 \begin{eqnarray*}
-Ric_{x_0}(v, v)=H''+{H'^2\over n}\leq -K|v|^2.
\end{eqnarray*}
That is to say $Ric\geq K$. 
Moreover, we 
can conclude that the inequality $(\ref{hhhh2})$ becomes an equality for all geodesics linking all $\rho_0$ and $\rho_1$ in 
the Wasserstein space $P_2(M, g, \nu)$  if and only if $(M, g)$ is Einstein, i.e., $Ric=Kg$.

In the case $m>n$, we can use the same argument as above to prove that the inequality $(\ref{hhhh0H})$ or $(\ref{hhhh0N})$  holds if and only if the CD$(K, m)$-condition holds on $(M, g, \mu)$, i.e., $Ric_{m, n}(L)\geq Kg$. 
Moreover, we can prove that the inequality $(\ref{hhhh0H})$ becomes an equality for  
all geodesics linking all $\rho_0$ and $\rho_1$ in 
the Wasserstein space $P_2(M, g, \nu)$ 
if and only if $(M, g)$ is $(K, m)$-Einstein, i.e., $Ric_{m, n}(L)=Kg$.

\medskip

In general case of $p\neq 1$, recall that
\begin{eqnarray}
H_p''+{1\over m}H_p'^2
&=& \left(  p-1+{1\over m} \right) \left[\left(  \int_M L \phi d\gamma \right)^2-   \int_M  |L \phi|^2 d\gamma \right]-\int_M Ric_{m, n}(L)(\nabla \phi, \nabla \phi) d\gamma\nonumber\\
& &\hskip1cm -\int_M \left[ {m-n\over mn}\left(L \phi+{m\over m-n}\nabla V\cdot\nabla \phi\right)^2 + \left\|\nabla^2 \phi-{\Delta \phi\over n} g\right\|^2_{\rm HS} \right]d\gamma.\nonumber\\
& &\label{NKsigma1}
\end{eqnarray}

In particular, when $m=n$, 
\begin{eqnarray}
H_p''+{1\over n} H_p'^2
&=& \left(  p-1+{1\over n} \right) \left[\left(  \int_M \Delta \phi d\gamma \right)^2-   \int_M  |\Delta\phi|^2 d\gamma \right]-\int_M Ric(\nabla \phi, \nabla \phi) d\gamma\nonumber\\
& &\hskip1cm -\int_M  \left\|\nabla^2 \phi-{\Delta \phi\over n} g\right\|^2_{\rm HS} d\gamma.\nonumber\\
& &\label{NKsigma2}
\end{eqnarray}

This yields, if $Ric\geq K$, we have 
\begin{eqnarray}
H_p''+{1\over n} H_p'^2
\leq -K \int_M |\nabla \phi|^2d\gamma.\label{NKsigma3}
\end{eqnarray}
By the fact 
$${d^2\over dt^2}N_{n, p}(\rho(t))=\left({1\over n} H_p''+{1\over n^2}H_p'^2\right)N_{n, p}(\rho(t)),$$
we see that $(\ref{NKsigma3})$ is equivalent to 
\begin{eqnarray}
{d^2\over dt^2}N_{n, p}(\rho(t))
\leq \left(-{1\over n} K \int_M |\nabla \phi|^2d\gamma\right)N_{n, p}(\rho(t)).\label{NKsigma4}
\end{eqnarray}

Conversely, we can prove that  if $(\ref{NKsigma3})$ or  $(\ref{NKsigma4})$ holds, then
\begin{eqnarray*}
Ric\geq Kg.
\end{eqnarray*}
For this, notice that $(\ref{NKsigma2})$ holds at $t=0$ for the solution $(\rho(t), \phi(t))$ of the continuity equation and the Hamilton-Jacobi equation with any initial data 
$(\rho(0), \phi(0))$. Taking a sequence of $d\gamma_k(0)={\rho_k(0)^pd\nu\over \int_M \rho_k(0)^pd\nu}$ converges weakly to the Dirac mass at $x_0\in M$, we derive that at $t=0$

\begin{eqnarray*}
H_p''+{1\over n} H_p'^2=-\left[Ric_{x_0}(v, v)+\left\|\nabla^2 \phi(x_0)- {\Delta \phi(x_0)\over n} g\right\|^2_{\rm HS}\right].
\end{eqnarray*}
For any fixed $v\in T_{x_0}M$, choose  $\phi$ such that 
\begin{eqnarray*}
\nabla \phi(x_0)&=&v,\\
\nabla^2\phi(x_0)&=&{\Delta \phi(x_0)\over n}I_n.
\end{eqnarray*}
This implies that

\begin{eqnarray*}
H_p''+{1\over n} H_p'^2=-Ric_{x_0}(v, v).
\end{eqnarray*}
Therefore, if $(\ref{NKsigma3})$ or $(\ref{NKsigma4})$ holds, we have  

 \begin{eqnarray*}
Ric_{x_0}(v, v)=-H_p''-{1\over n} H_p'^2\geq K|v|^2.
\end{eqnarray*}
That is to say $$Ric\geq Kg.$$ 

Moreover, in the case where $Ric=Kg$, we 
can conclude that the equality holds in $(\ref{NKsigma3})$, i.e., 
\begin{eqnarray}
H_p''+{1\over n} H_p'^2
=-K \int_M |\nabla \phi|^2d\gamma,\label{NKsigma3b}
\end{eqnarray}
or equivalently, the equality holds in $(\ref{NKsigma4})$,  i.e., 
\begin{eqnarray}
{d^2\over dt^2}N_{n, p}(\rho(t))
= \left(-{K\over n}  \int_M |\nabla \phi|^2d\gamma\right)N_{n, p}(\rho(t))\label{NKsigma4b}
\end{eqnarray}
holds if and only if $\phi$ is 
a Hessian soliton 
\begin{eqnarray*}
\nabla^2\phi={ \Delta \phi \over n} g.
\end{eqnarray*}

Similarly, we can prove that the inequality 
\begin{eqnarray}
H_p''+{1\over m} H_p'^2
\leq -K \int_M |\nabla \phi|^2d\gamma,\label{NKsigma3b}
\end{eqnarray}
equivalently, 
\begin{eqnarray}
{d^2\over dt^2}N_{m, p}(\rho(t))
\leq \left(-{K\over m} \int_M |\nabla \phi|^2d\gamma\right)N_{m, p}(\rho(t))\label{NKsigma4c}
\end{eqnarray}
holds for any $\rho_0$ and $\rho_1$ in the Wasserstein space $P_2(M, g, \mu)$
if and only if the CD$(K, m)$-condition holds on $(M, g, \mu)$, i.e., $Ric_{m, n}(L)\geq Kg$. 
Moreover, in the case where $Ric_{m, n}(L)=Kg$, we can prove that the equality holds in 
$(\ref{NKsigma3b})$ (equivalently, $(\ref{NKsigma4c})$) 
if and only if $\phi$ is a Hessian soliton 
\begin{eqnarray*}
\nabla^2\phi={ \Delta \phi \over n} g,
\end{eqnarray*}
and
\begin{eqnarray*}
 {m-n\over mn}\left(L \phi+{m\over m-n}\nabla V\cdot\nabla \phi\right)^2=0.
\end{eqnarray*}

We have therefore finished the proof of the rigidity part of Theorem \ref{Thm2}. \hfill $\square$

\subsection{Summary and further problems}

Theorem \ref{Thm1} and Theorem \ref{Thm2} indicate the equivalence between the curvature-dimension CD$(K, m)$-condition, entropy differential inequalities (EDI) and entropy power differential  inequalities (EPDI) 
along the geodesics on the Wasserstein space over Riemannian manifolds. In particular, each of the equivalent conditions (ii)-(vii) in Theorem \ref{Thm2} is equivalent to the definition inequality $(\ref{S(K, N)})$ (for $N'\geq N=m$) 
which was used by Sturm \cite{Sturm06b} to define the CD$(K, N)$-condition on metric measure spaces. They are also equivalent to Lott-Villani's characterization of CD$(K, N)$ condition using functionals in  
the class $\mathcal{DC}_N$. Therefore, we can 
use each of (ii)-(vii) or its integral forms to introduce definition of CD$(K, N)$-condition on 
metric measure spaces. In particular, all these equivalent conditions are stable under the Gromov-Hausdorff convergence on metric measure spaces.

Theorem \ref{rigidity3} and Theorem \ref{rigidity5} also indicate that the rigidity models for the enhanced entropy 
differential inequality and the enhanced entropy power differential inequality 
along all smooth 
geodesics on the Wasserstein space over Riemannian manifolds are $K$-Einstein manifolds 
or $(K, m)$-Einstein manifolds.  The proofs of rigidity theorems need to use the Bochner formula 
 $(\ref{BWF})$ rather than the Bochner inequality $(\ref{U2c})$.

Theorem \ref{W-entropy}  extends the monotonicity 
and rigidity theorems of the $W$-entropy 
associated with the Shannon entropy \cite{LL18, LL24} to the R\'enyi entropy for the geodesics on the Wasserstein space over complete 
Riemannian manifolds with CD$(0, m)$-condition.  

In \cite{EKS},  Erbar-Kuwada-Sturm  proved that the Shannon entropy power $
\exp\left(-{2\over N}{\rm Ent}(\rho))\right)$ is concave along the heat equation on RCD$(0, N)$ metric measure spaces.  In \cite{KuwadaLi}, K. Kuwada and the author of this paper proved the monotonicity  of the $W$-entropy for heat equation associated with the Witten Laplacian on RCD$(0, N)$ metric measure spaces. Moreover, a rigidity theorem in a weak form was also proved for the $W$-entropy for heat equation on RCD$(0, N)$ metric measure spaces.  

\medskip

It is interesting to  raise some problems for further study in the future.

\begin{problem}
 (i)  Can we extend 
 Theorem \ref{Thm1}, Theorem \ref{Thm2} and Theorem \ref{W-entropy} (in particular, the rigidity theorems) to geodesics on the Wasserstein space over RCD$(0, N)$ or RCD$(K, N)$ metric measure spaces?
 
(ii) Can we extend 
 Theorem \ref{Thm1}, Theorem \ref{Thm2} and Theorem \ref{W-entropy} 
(in particular, the rigidity theorems)  
 to geodesics on the Wasserstein space over RCD metric measure spaces with time-dependent metrics and measures?

(iii) How to introduce the notion of $n$-dimensional metric measure spaces with 
constant $N$-Ricci curvature, i.e., $Ric_{n, N}=K$?

(iv) Can we extend EDI, EPDI and the $W$-entropy formula to geodesics on the Wasserstein space over metric measure spaces with $(K, N)$-super Ricci flows on non-smooth metric measure spaces introduced by 
Kopfer and Sturm  \cite{Sturm18a} and Sturem \cite{Sturm18b} ?

(v) How to introduce the notion of the $(K, n, N)$-Ricci flow of $n$-dimensional RCD metric measure spaces with ${1\over 2}{\partial g\over \partial t}+Ric_{n, N}=Kg$?

\end{problem}

See also \cite{Li25WRF, LiZhang}. 

\medskip

\noindent{\bf Acknowledgement}. The author would like  to thank Prof. N. Mok for suggestion,  Prof. B-X. Han and 
Dr. S. Li for helpful discussions, Dr. R. Lei and Dr. Y.-Z. Wang for careful reading on the manuscript. 
He is also grateful to an anonymous referee for pointing out an mistake in the formulation of the rigidity theorem in 
the first version of this paper.

\begin{flushleft}
\medskip\noindent

Xiang-Dong Li, State Key Laboratory of Mathematical Sciences, Academy of Mathematics and Systems Science, Chinese
Academy of Sciences, No. 55, Zhongguancun East Road, Beijing, 100190,  China\\
E-mail: xdli@amt.ac.cn

and

School of Mathematical Sciences, University of Chinese Academy of Sciences, Beijing, 100049, China
\end{flushleft}


\begin{thebibliography}{99}
\bibitem{AG} L. Ambrosio, N. Gigli, A User's Guide to Optimal Transport Theory, 
CIME Lecture Notes in Mathematics, Springer, 2012.
\bibitem{BE} D. Bakry, M. Emery, Diffusion hypercontractives, S\'em. Prob. XIX, Lect. Notes in Maths. 1123 (1985), 177-206.
\bibitem{BGL} D. Bakry, Y. Gentil, M. Ledoux, Analyse and Geometry for Markov Diffusion Operators, Springer, 2014. 
\bibitem{Bla} N.M. Blachman, The convolution inequality for entropy powers. IEEE Trans. Inform. Theory 2, 267-271, (1965).
\bibitem{BB}  J.-D. Benamou, Y. Brenier, A computational fluid mechanics solution to the Monge-Kantorovich mass
transfer problem. Numer. Math. 84, 375-393 (2000)
\bibitem{CZ} H.-D. Cao, X.-P. Zhu, A complete proof of the Poincar\'e and geometrization conjectures: application of the Hamilton--Perelman theory of the Ricci flow, Asian J. Math. 10 (2006), no. 2, 165-492, and Asian J. Math. 10 (2006), no. 4, 663.
\bibitem{CLN} B. Chow, P. Lu, L. Ni, Hamilton's Ricci flow, Lectures in Contemporary Mathematics, 3, Amer. Math. Soc. Providence, Rhode Island,  Beijing, 2006.
\bibitem{CCLLN} B. Chow, S.-C. Chu,  D. Glickenstein,  C. Guenther,  J. Isenberg,  T. Ivey, D. Knopf,  P. Lu,  F. Luo and L. Ni, The Ricci flow: Techniques and Applications,  Part I, Geometric aspects, Mathematical Surveys and Monographs, 135. American Mathematical Society, Providence, RI, 2007.
\bibitem{CEMS} D. Cordero-Erausquin, R.J.  McCann,  and M.  Schmuckenschl\"ager, 
A Riemannian interpolation inequality \`a la Borell, Brascamp and Lieb. Invent. Math. 146, 2 (2001), 219-257.   
\bibitem{Cost} M.H.M. Costa, A new entropy power inequality, IEEE Trans. Inf. Theory, IT- 31, (6) 751-760, (1985).
\bibitem{CT} T. M. Cover, J.A. Thomas, Elements of Information Theory, Wiley-Interscinces, A John Wiely and Sons, INC., Publication. 
\bibitem{Dem1} A. Dembo, A simple proof of the concavity of the entropy power with respect to the variance of additive normal noise, IEEE Trans. Inform. Theory 35, 887-888, (1989).
\bibitem{Dem2}  A. Dembo, T.M. Cover, and J.A. Thomas, Information theoretic inequalities, IEEE Trans. Inf. Theory, 37, (6), 1501-1518, (1991).
\bibitem{EKS} M. Erbar, K. Kuwada, K.-T. Sturm, On the equivalence of the
  entropic curvature-dimension condition and Bochner's inequality on metric
  measure spaces, Invent.\ Math. 201 (2015), no.~3, 993--1071.
  \bibitem{FLZ} F. Fang, X.-D. Li and Z. Zhang, Two generalizations of the Cheeger-Gromoll splitting theorem via the
Bakry-Emery Ricci curvature. Ann. Inst. Fourier 59(2) (2009), 563--573.
\bibitem{Ka} T. Kato, The Cauchy problem for quasi-linear symmetric hyperbolic systems, Arch Rational Mech Anal. 58 (1975), 181--205. 
\bibitem{GW}  R. E. Greene and H. Wu, $C^\infty$-convex functions and manifolds of positive curvature, 
  Acta Math. 137 (1976), no. 3-4, 209--245. 
\bibitem{H82} R.S. Hamilton, Three manifolds with positive Ricci curvature. J. Diff. Geom. 17 (1982), 255-306.
\bibitem{H95}R.S. Hamilton, The formation of singularities in the Ricci flow, Surveys in Differential
Geometry, 2, 7-136, International Press, 1995.
\bibitem{KL} B. Kleiner, J. Lott, Notes on Perelman's papers, Geom. Topol. 12 (2008),no. 5, 2587-2855.
\bibitem{Sturm18a}  E. Kopfer, K.-Th. Sturm, Heat flow on time-dependent metric measure spaces and super-Ricci flows. Comm. Pure Appl. Math. 71 (2018), no. 12, 2500-2608.
\bibitem{KuwadaLi}  K. Kuwada, X.-D. Li, Monotonicity and rigidity of the W-entropy on RCD(0;N)
spaces, Manuscripta Math. 164 (2021), 119-149.
\bibitem{LL15}S. Li, X.-D. Li,  $W$-entropy formula for the Witten Laplacian on
manifolds with time dependent metrics and potentials,  Pacific J. Math. Vol. 278 (2015), No. 1, 173-199.
\bibitem{LL18a} S. Li and X.-D. Li, Hamilton differential Harnack inequality and $W$-entropy for Witten Laplacian on Riemannian manifolds, J. Funct. Anal. 274 (2018), 3263-3290.
\bibitem{LL18} S. Li, X.-D. Li, $W$-entropy formulas on super Ricci flows and Langevin deformation on Wasserstein space over Riemannian
manifolds. Sci China Math, 2018, 61: 1385-1406.
\bibitem{LL19} S. Li and X.-D. Li, $W$-entropy, super Perelman Ricci flows and $(K, m)$-Ricci solitons,  J. Geom. Anal. 30 (2020), no.~3, 3149--3180.
\bibitem{LL24}   S. Li, X.-D. Li,  $W$-entropy formulas and Langevin deformation of flows on Wasserstein space over Riemannian manifolds. Probab. Theory and Related Fields, 188 (2024), 911-955.
\bibitem{LL20a} S. Li, X.-D. Li, On the Shannon entropy power on Riemannian manifolds and Ricci flows, Tohoku Math. J. (2) 76 (2024), 577-608. 
\bibitem{Li12} X.-D. Li, Perelman’s entropy formula for the Witten Laplacian on Riemannian manifolds via Bakry-Emery Ricci curvature, Math. Ann. 353 (2012), 403–437.
\bibitem{Li25WRF}X.-D. Li,, On Perelman’s W-entropy and Shannon entropy power for super Ricci flows on metric measure spaces, arXiv:2505.03202v1
\bibitem{Li26} X.-D. Li, On the R\'enyi entropy power and the $W$-entropy for nonlinear diffusion equations on Riemannian manifolds, preprint, 2026.
\bibitem{LiZhang} X.-D. Li, E. Zhang,  On the W -entropy and Shannon entropy power on 
RCD$(K, N)$ and RCD$(K, n, N)$ spaces, arxiv2504.01864
\bibitem{Lo1} J. Lott, Some geometric calculations on Wasserstein space, Comm. Math. Phys. 277 (2008), no. 2, 423-437.
\bibitem{Lo2} J. Lott, Optimal transport and Perelman's reduced volume, Calc. Var. and Partial Differential Equations 36, 49-84 (2009).
\bibitem{LV} J. Lott, C. Villani, Ricci curvature for metric-measure spaces via optimal transport, Annals of Math. 169, 903-991, 2009.
\bibitem{Ot} F. Otto, The geometry of dissipative evolution equations: the porous medium equation, Commun. in Parial Differential Equations 26 (1 and 2), 101-174 (2001).
\bibitem{Mc} R. McCann, Polar factorization of maps on Riemannian manifolds, Geometric and Functional Analysis, 11, 3(2001), 589-608.
\bibitem{MT} J. W. Morgan, G. Tian, Ricci flow and the Poincar\'e
conjecture, Clay Mathematics Monographs, 3. American Mathematical
Society, Providence, RI; Clay Mathematics Institute, Cambridge, MA,
2007.
\bibitem{N1} L. Ni, The entropy formula for linear equation, J.
Geom. Anal. 14 (1), 87-100, (2004).
\bibitem{N2} L. Ni, Addenda to ``The entropy formula for linear
equation'', J. Geom. Anal. 14 (2), 329-334, (2004).

\bibitem{P1} G. Perelman, The entropy formula for the Ricci flow and its geometric applications, http://arXiv.org/abs/maths0211159.
\bibitem{ST}G. Savar\'e, G. Toscani,  The concavity of R\'enyi entropy power, IEEE Transactions on Information Theory Vol 60, Issue: 5, May 2014, 2687-2693. DOI: 10.1109/TIT.2014.2309341
\bibitem{Sh48}  C. E. Shannon,  A mathematical theory of communication, Bell Syst. Tech. J., vol. 27, pp. 623-656, Oct. 1948.
\bibitem{Stam}  A. J. Stam, Some inequalities satisfied by the quantities of information of Fisher and Shannon, Information and Control, vol. 2, pp. 101-112, June 1959.
\bibitem{Sturm05}  K. T. Sturm, Convex functionals of probability measures and nonlinear diffusions on manifolds,  J. Math. Pures Appl. 84 (2005) 149-168.
\bibitem{Sturm06a} K.T. Sturm, On the geometry of metric measure spaces. {\rm I}, Acta.\ Math. \textbf{196} (2006), no.~1, 65--132.
  \bibitem{Sturm06b} K.T. Sturm, On the geometry of metric measure spaces. {\rm II}, Acta.\ Math. \textbf{196} (2006), no.~1, 133--177.
\bibitem{Sturm18b} K.-Th. Sturm, Super-Ricci flows for metric measure spaces. J. Funct. Anal. 275 (2018), no. 12, 3504–3569.
\bibitem{RS} K.-T. Sturm, M.-K. von Renesse, Transport inequalities, gradient estimates, entropy and Ricci curvature. Commun. Pure Appl. Math. 58(7), 923-940 (2005)
\bibitem{V0} C. Villani, A Short Proof of the Concavity of Entropy Power,  IEEE Transactions on Information Theory, Vol. 46, No. 4,1695-1696, July 2000.
\bibitem{V1} C. Villani,  Topics in Mass Transportation, Grad. Stud. Math., Amer. Math. Soc., Providence,
RI, 2003.
\bibitem{V2}C. Villani, Optimal Transport, Old and New, Springer, 2008.
\end{thebibliography}
\end{document}